\documentclass [11pt]{amsart}
\usepackage{amsmath,mathrsfs,amssymb}
\usepackage{ulem}
\usepackage{tikz}
\usepackage{color}
\usepackage{pstricks}
\usepackage{graphicx}
\usepackage{epstopdf}
\usepackage{young}
\usepackage[vcentermath]{youngtab}
\newtheorem{Nota}{Remark}

\title{Orthogonal polynomials associated with the deltoid curve }
\author{Olfa Zribi }
\address{Institut de Math\'ematiques de Toulouse, Universit\'e Paul Sabatier, 118, route de Narbonne, Toulouse, France\\ D\'epartement de Math\'ematiques, Universit\'e El Manar, Tunis}
\begin{document}
\newcommand{\un}{{\bf 1}}
\newcommand{\va}{variable al\'eatoire}
\newcommand{\vas}{variables al\'eatoires}
\newcommand{\R}{{\mathbb {R}}}
\newcommand{\Q}{{\mathbb {Q}}}
\newcommand{\N}{{\mathbb {N}}}
\newcommand{\Z}{\mathbb {Z}}
\newcommand{\pr}{{\mathbb P}}
\newcommand{\Nbar}{\overline{\mathbb {N}}} 
\newcommand{\ds}{\frac{d}{ds}}
\newcommand{\fp}{{\mathcal F}_p}
\newcommand{\su}{\sum_{i=1}^{d}}
\newcommand{\sij}{\sum_{i,j=1}^{d}}
\newcommand{\dx}{\Delta_{p+1}^h\overline{x}}
\newcommand{\eE}{\stackrel{E}{=}}
\newcommand{\tx}{\tilde{x}}
\newcommand{\tda}{{\it temps d'arr\^et}}
\newcommand{\wee}{\wedge }
\newcommand{\Tr}{{\mathcal T} }
\newcommand{\TR}{{\mathcal T R} }
\newcommand{\Rec}{{\mathcal R} }
\newcommand{\tP}{{}^{t}P }
\newcommand{\bex}{\begin{exemp}}
\newcommand{\eex}{\end{exemp}}
\newcommand{\bexs}{\begin{exemps}}
\newcommand{\eexs}{\end{exemps}}
\newcommand{\bprop}{\begin{proposition}}
\newcommand{\eprop}{\end{proposition}}
\newcommand{\bthm}{\begin{theoreme}}
\newcommand{\ethm}{\end{theoreme}}
\newcommand{\bcor}{\begin{corollaire}}
\newcommand{\ecor}{\end{corollaire}}
\newcommand{\blem}{\begin{lemme}}
\newcommand{\elem}{\end{lemme}}
\newcommand{\beq}{\begin{equation}}
\newcommand{\eeq}{\end{equation}}
\newcommand{\beqna}{\begin{eqnarray}}
\newcommand{\eeqna}{\end{eqnarray}}
\newcommand{\beqnas}{\begin{eqnarray*}}
\newcommand{\eeqnas}{\end{eqnarray*}}

\newcommand{\bmat}{\begin{pmatrix}}
\newcommand{\emat}{\end{pmatrix}}

\newcommand{\cL}{{\mathbb {L}}}

\newcommand{\bpf}{\begin{preuve}}
\newcommand{\epf}{ \end{preuve} \medskip}

\newcommand{\benum}{\begin{enumerate}}
\newcommand{\eenum}{\end{enumerate}}

\newcommand{\bitem}{\begin{itemize}}
\newcommand{\eitem}{\end{itemize}}

\newcommand{\brmq}{\begin{rmq}}
\newcommand{\ermq}{\end{rmq}}

\newcommand{\brmqs}{\begin{rmqs}}
\newcommand{\ermqs}{\end{rmqs}}

\newcommand{\bapp}{\begin{application}}
\newcommand{\eapp}{\end{application}}

\newcommand{\bapps}{\begin{applications}}
\newcommand{\eapps}{\end{applications}}

\newcommand{\bdefi}{\begin{definition}}
\newcommand{\edefi}{\end{definition}}

\newcommand{\bequa}{\begin{equation}}
\newcommand{\eequa}{\end{equation}}

\newcommand{\dR}{\ensuremath{\mathbb{R}}}
\newcommand{\p}[4]{{#3}\!\left#1{#4}\right#2}

\newcommand{\proofend}{~$\blacktriangleright$}
\newcommand{\proofbegin}{~$\blacktriangleleft$}
\newcommand{\ABS}[1]{\ensuremath{{\left| #1 \right|}}} 
\newcommand{\BRA}[1]{\ensuremath{{\left\{#1\right\}}}} 
\newcommand{\DP}[1]{\ensuremath{{\left<#1\right>}}} 
\newcommand{\NRM}[1]{\ensuremath{{\left\Vert #1\right\Vert}}} 
\newcommand{\SBRA}[1]{\ensuremath{{\left[#1\right]}}} 
\newcommand{\varf}[1]{\mathbf{Var}_{#1}}
\newcommand{\var}[2]{\p(){\varf{#1}}{#2}}
\newcommand{\disp}{\displaystyle}
\newcommand{\e}{\varepsilon}

\newcommand{\bentf}[1]{\mathbf{H}_{#1}}
\newcommand{\bent}[2]{\p(){\bentf{#1}}{#2}}
\newcommand{\eentf}{\mathbf{N}}
\newcommand{\eent}[1]{\p(){\eentf}{#1}}
\newcommand{\entf}[1]{\mathbf{Ent}_{#1}}
\newcommand{\ent}[2]{\p(){\entf{#1}}{#2}}
\newcommand{\ientf}[1]{\mathbf{H}^{#1}}
\newcommand{\ient}[2]{\p(){\ientf{#1}}{#2}}
\newcommand{\PAR}[1]{\ensuremath{{\left(#1\right)}}} 
\newcommand{\CC}{\ensuremath{\mathcal C }}
\newcommand{\UU}{\ensuremath{\mathcal U }}
\newcommand{\A}{{\mathcal A}}

\newcommand{\PT}[1]{\mathbf{P}_{\!#1}}
\newcommand{\Pt}[1][t]{\ensuremath{\mathbf{P}_{\!#1}}}
\newcommand{\NPsi}{\rm{N}\Psi }

\renewcommand{\red}{\color{red}\tt}
\renewcommand{\blue}{\color{blue}}
\renewcommand{\green}{\color{green}}
\newcommand{\oran}{\color{orange}}

\definecolor{green}{rgb}{0,.7,.2}
\definecolor{orange}{rgb}{0.9,.5,0}
\newcommand{\bx}{{\bf x}}
\newcommand{\by}{{\bf y}}
\newcommand{\bO}{{\bf 0}}
\newcommand{\bGam}{{\bf \Gamma}}
\newcommand{\bu}{{\bullet}}
\newcommand{\ind}{{\bf 1}}
\newcommand{\HH}{\mathbb {H}}
\newcommand{\hh}{{\rm H}}
\newcommand{\sph}{{\mathbb {S}}}
\newcommand{\LL}{{\rm L}}
\newcommand{\pt}{{(P_t)}_{t\geq 0}}
\newcommand{\Var}{{\rm Var}} 
\newcommand{\sob}{{\rm Sob}} 

\def\I{1\!{\rm l}}
\def\C{{\mathbb{C}}}
\def\E{{\mathbb{E}}}
\def\tr{\textmd{trace}\,}
\def\grad{\textmd{grad}\,}
\def\ssi{\Longleftrightarrow}
\def\lag{\langle}
\def\rag{\rangle}
\def\OU{\textsc{Ornstein-Uhlenbeck}}
\def\Lip{\textsc{Lipschitz}}
\def\logsob{logarithmic \textsc{Sobolev}}
\def\Sob{ \textsc{Sobolev}}
\def\Mkv{ \textsc{Markov}}
\def\LapBel{ \textsc{Laplace}-\textsc{Beltrami}}
\def\Riem{ \textsc{Riemann}}
\def\Ric{ {\rm{Ric}}} 
\def\det{{ \rm{det}}}  
\def\Ent {{\rm{Ent}}}   
\def\essup{ \rm{esssup}}
\def\Sob{ \rm{Sob}}
\def\SPhi{ \rm{S}\Phi}
\def\dis{\displaystyle }
\def\dd{\partial }
\def\Id{{\rm{Id}}} 
\def\div{{\rm{div}}}
\def\sc{{\rm{sc}}}  
\def\eps{\varepsilon} 
\def \capa {{\rm Cap}} 

\def\cA{{\mathcal A }}
\def\cB{{\mathcal B }}
\def\cC{{\mathcal C}}
\def\cD{{\mathcal D}}
\def\cE{{\mathcal E }}
\def\cF{{\mathcal  F}}
\def\cG{{\mathcal  G}}
\def\cH{{\mathcal  H}}
\def\cL{{\mathcal  L}}
\def\cS{{\mathcal  S}}
\def\cT{{\mathcal  T}}
\def\cK{{\mathcal  K}}
\def\cP{{\mathcal P }}
\def\cR{{\mathcal R }}
\def\cV{{\mathcal V}}
\def\cO{{\mathcal O}}

\def\bbR{\mathbb{R}}
\def\bbC{\mathbb{C}}
\def\bbZ{\mathbb{Z}}
\def\bbN{\mathbb{N}}
\def\bbQ{\mathbb{Q}}
\def\bbL{\mathbb{L}}

\def\bpm{\begin{pmatrix}}
\def\epm{\end{pmatrix}}
\def\page{\newpage}

\newtheorem{theoreme}{Theorem}[section]
\newtheorem{lemme}[theoreme]{Lemma}
\newtheorem{definition}[theoreme]{Definition}
\newtheorem{proposition}[theoreme]{Proposition}
\newtheorem{corollaire}[theoreme]{Corollary}
\newtheorem{example}[theoreme]{Example}
\newtheorem{fthm}{Theorm}[section]
\newtheorem{flem}[fthm]{Lemma}

\newenvironment{exemp}{\noindent{\bf Example. --- }}{\par}
\newenvironment{exemps}{\noindent{\bf Examples}\benum}{\eenum\par}
\newtheorem{rmq}[theoreme]{Remark}
\newtheorem{rmqs}[theoreme]{Remarks}

\newenvironment{preuve}{\noindent{\it Proof. --- }}
{\hfill\rule{1.3mm}{2mm}\par} 
\newenvironment{application}{\noindent{\bf Application. --- }}{\par}
\newenvironment{applications}{\noindent{\bf Applications. --- 
}\benum}{\eenum\par}

\theoremstyle{definition}
\newtheorem{fdefi}[fthm]{Definition}
\newtheorem{frem}[fthm]{Remark}
\newtheorem{fcor}[fthm]{Corollary}

\def\bZ{\overline{Z}}
\def\bz{\overline{z}}
\maketitle
\begin{abstract}
We study  a family of bivariate orthogonal polynomials associated  to the deltoid curve. These polynomials arise when classifying bivariate diffusion operators that have discrete spectral decomposition given by orthogonal polynomials with respect to some compactly-supported probability measure on the interior of the deltoid curve.
\end{abstract}



\section{Introduction}

Orthogonal polynomials in the interior of the deltoid curve is one example of the 11 families of orthogonal polynomials on a compact domain in dimension 2 which are at the same time eigenvectors of an elliptic diffusion operators, see \cite{BOZ}. It is also one of the most intriguing one, and have been put forward by Koornwinder \cite{Koorn1, Koorn2, Koorn}, (see also \cite{Dunkel}, \cite{KrallS}).  This family of bivariate polynomials depends on a parameter $\alpha$ : if $P(X,Y)$ is the algebraic equation of the boundary of $\Omega$,  so that $P(X,Y)>0$ on $\Omega$, then the measure has density $C_\alpha P(X,Y)^\alpha$ with respect to the Lebesgue measure. Here, $\alpha > -5/6$, as we shall see below (Proposition~\ref{prop.int.mesure}).

The two special cases $\alpha=1/2$ and $\alpha=-1/2$ play a special role   and have been particularly investigated, see \cite{beerends, dunn} and also  \cite{Pinsky1,Pinsky} for a spectral point of view. The first one is the image of the Euclidean Laplace operator through the symmetries of the triangular lattice, and the second one is the image of the Casimir operator on $SU(3)$ through  the spectral decomposition.  It is of course not a surprise, since the root system of $SU(3)$ is $A_2$, which corresponds to the triangular lattice (see {\cite{Bourb.Lie1, Bourb.Lie2, Bourb.lie3, yvette,hum, opdam}}).  Those two cases are referred below as the geometric cases.
The analysis of the geometric cases provides  some insight on the general model. Therefore, although this aspect is quite well documented (see  \cite{beerends}, \cite{dunn}), we present  them in detail from the point of view of symmetric diffusion operators
 for the sake of completeness. Moreover, these models will provide us efficient insights towards the general situation, since  the study of
the general family is more delicate. In the core of the paper, we derive recurrence formulae for the generic measures, which, as is the case  in dimension 1, turn out to be a 3 term recurrence formula (although such a simple form is not to be expected in general in dimension 2). Those recurrence formulae  take a particularly  simple form in the two geometric cases.

As mentioned above, the study of the  two specific geometric cases  lead to simple representations of the eigenvectors. In the Euclidean case ($\alpha = -1/2$), we  have a very simple presentation of the eigenvectors. The $SU(3)$ case  ($\alpha = 1/2$) leads to a representation of the eigenvectors through the representation of the symmetric group, associated with Young diagrams (see \cite{James}). Finally, we  derive in the general case  partial generating functions, leading to another representation of the orthogonal polynomials, which also provides a complete generating function in the two geometric cases.

The paper is organized as follows : Section~\ref{sec1} is a short presentation of the general setting of symmetric diffusion processes associated with orthogonal polynomials, mostly inspired from \cite{bglbook} and \cite{BOZ}. In Section~\ref{sec.deltoid}, we give the explicit formulae for the measure and the generator associated to the Deltoid model, and introduce the complex variables in which the operator  takes a much simpler form, leading to the explicit values for the eigenvalues of the operator. Section~\ref{sec.laplacien2d} is the presentation of the Euclidean case (that is the case $\alpha= -1/2$), while Section~\ref{sec.su3} presents the $SU(3)$ case, that is $\alpha= 1/2$. Section~\ref{sec.eigen} concentrates on   recurrence formulae in the general case, and Section~\ref{sec.other} provides some representation of the eigenvectors, first  in the geometric cases and then (but only partially for the generating functions) in the general case.

\section{Orthogonal polynomials and diffusion generators\label{sec1}}
Let $\Omega$ be an open  bounded domain in $\mathbb{R}^d, d \geq 1$, with piecewise smooth boundary, and let $\mu$  a probability measure on $\overline{\Omega}$. Recall from p. 32 in \cite{Dunkel} that a family of polynomials $P_{\tau}: \mathbb{R}^d \mapsto \mathbb{R}$ is orthogonal in $L^2(\Omega, \mu)$ if
\begin{equation*}
\int P_{\tau}(x)P_{\tau'}(x)\mu(dx) = 0
\end{equation*}
where $\tau = (\tau_1,\cdots,\tau_d) \in \mathbb{N}^d$ is a multi-index, whenever $|\tau| : = \tau_1+\dots + \tau_d \neq \tau'_1+\dots \tau'_d = |\tau'|$. In contrast to the real one variable setting, this family needs not to be unique in higher dimensions due to various orders one may choose when applying the Gram-Schmidt process to the canonical basis $(x_1^{\tau_1}\dots x_d^{\tau_d})_{\tau \in \mathbb{N}^d}$ (see the bottom of p.31 in \cite{Dunkel}). However, in many situations, there are natural choices for this family orthogonal polynomials. In particular, it may happen that they are also eigenvectors of some diffusion differential operator. This is the case for the classical family of orthogonal polynomials in dimension 1, Hermite, Laguerre and Jacobi (although only the last one corresponds to a bounded domain, see \cite{Mazet97}).

On the other hand,   when solving stochastic differential equations in probability theory, one is often led to consider second order differential operators on $\Omega$ which are symmetric with respect in $\cL^2(\mu)$, at least when one restricts it's attention to the set $\cC^\infty_c(\Omega)$ of smooth functions compactly supported in $\Omega$.  When $\mu$ has a smooth positive density on $\Omega$, these operators may be represented  as
\beq\label{eqdiffusion1}
\LL f := \frac{1}{\rho} \sum_{k,j=1}^d \partial_k\left(g_{kj}\rho\partial_j f\right)
=  \sum_{k,j=1}^d g_{kj}\partial_{kj}^2 f + \sum_{j=1}^d b_j\partial_j f
\eeq
where $g = (g_{kj}(x))_{k,j=1}^d, x \in \Omega$ is  a symmetric non negative  matrix depending smoothly on $x\in \Omega$ and
\begin{equation*}
b_j = \frac{1}{\rho} \sum_{k=1}^d\partial_k(g_{kj}\rho), \quad j \in \{1,\dots,d\}.
\end{equation*}
The coefficients $b_j(x)$ are called the drift terms of the operator $\LL$.

We call such operators symmetric diffusion operators. They are related  to Markov diffusion processes $(X_t)$ with values in $\Omega$ through the fact that for any smooth function $f$, the processes $f(X_t)-\int_0^t \LL f(X_s) \, ds$ is a (local) martingale. When the operator $\LL$ is essentially self-adjoint, then this entirely characterizes the law of the processus $(X_t)$ (at least as long as we only consider  the finite dimensional marginals). The operator $\LL$ is called the infinitesimal generator of the process $(X_t)$.

Working with such diffusion operators,   it is often convenient to introduce the so called carr\'e du champ operator
$$\Gamma(f,g)= \frac{1}{2}\big(\LL(fg)-f\LL g-g\LL f),$$ and observe that $\LL$ is entirely determined from the knowledge of $\Gamma$ and $\mu$ through the integration by parts formula
$$\int_\Omega f \LL g\, d\mu = \int_\Omega g\LL f\, d\mu = -\int_\Omega \Gamma(f,g)\, d\mu,$$ valid at least when $f$ and $g$ are smooth and compactly supported in $\Omega$. Moreover From the representation~\eqref{eqdiffusion1}, it is immediate that $b_i(x)= \LL(x_i)$ and $g_{ij}= \Gamma(x_i,x_j)$.

More generally, the change of variable formula, valid for any smooth $\Phi : \R^k\mapsto\R$, and any $k$-uple $f= (f_1, \cdots, f_k)$ of smooth functions, reads
\beq\label{change.var}\LL(\Phi (f)) = \sum_i \partial_i \Phi(f)\LL f_i + \sum_{ij} \partial^2_{ij} \Phi(f) \Gamma(f_i,f_j).\eeq

In particular, whenever for $i= 1, \cdots, k$, there exist functions $B_i$ and $G_{ij}$ such that $\LL f_i = B_i(f)$ and $\Gamma(f_i,f_j)= G_{ij}(f)$, one has
\beq
\LL \big(\Phi(f)\big) = (\LL_1 \Phi)(f)
\eeq
 where $\LL_1 $ is the new diffusion operator acting on the image of $\Omega$ under the function $f$ (which is not necessary a local diffeomorphism), as
$$ \LL_1 (\Phi) =\sum_{ij} G_{ij}(x)\partial^2_{ij} \Phi + \sum_i B_i(x)\partial_i \Phi,$$ which is called the image of $\LL$ under the function $f$.

 In the probabilistic interpretation, if $(X_t)$ is the stochastic process with generator $\LL$, then $f(X_t)$ is again a diffusion Markov process with generator $\LL_1$. In particular, the operator $\LL_1$ is symmetric with respect to the image measure of $\rho$ under the map $f$, which, following equation~\eqref{eqdiffusion1} may be often a efficient way to determine the image measure. When such situation occurs, we shall say that $\LL$ projects onto $\LL_1$.

In what follows, we restrict for simplicity to the case where the matrix $g(x)$ is positive definite on $\Omega$.  It is then natural to raise the question of determining when such $\LL$ may be extended as a self adjoint operator (see \cite{Yoshida}) with spectral decomposition given by a family of orthogonal polynomials with respect to $\mu$. In other words, one wants to determine for which choice of $\rho$ and $g$ there is a complete family of $\mu$-orthogonal polynomials which are at the same time eigenvectors for $\LL$. This will produce a natural choice for a basis of orthogonal polynomials.

It turns out that the general answer to this question is the following :

the functions $g_{ij}(x)$ are polynomials with degree at most two, and the boundary $\partial \Omega$ is included in the algebraic set $\{\det(g)=0\}$. More precisely if $P(x)$ denotes the irreducible equation of the boundary $\partial\Omega$, then there exists a family of degree $1$ polynomials $L_i(x)$ such that  for any $i$, the algebraic equation
\beq\label{eq.g.boundary}\sum_j g_{ij}\partial_j P= L_i P.\eeq Moreover, the sets of admissible density measures $\rho$  are entirely described from the algebraic structure of the boundary. In particular, when the determinant $\det(g)$ is irreducible, the only admissible density measures $\rho$ are $C(\lambda)\det(g)^\lambda$, for any real $\lambda$ such that $\rho^\lambda$ is $\cL^1(\Omega, dx)$, (see \cite{BOZ}). Once the boundary $\partial\Omega $ is given through it's irreducible equation, the coefficients $(g_{ij})(x)$ are entirely determined   from equation \eqref{eq.g.boundary}. It turns out that they are in general unique up to some scaling factor.

\section{The deltoid model\label{sec.deltoid}}
In dimension 2, up to affine transformations, there are only 11 bounded  sets $\Omega$ on which there exist a symmetric diffusion operator for which the associated eigenvectors are orthogonal polynomials with respect to the reversible measure (see { \cite{BOZ}) . One of the most intriguing one  is the interior of the deltoid curve, which is a  degree 4  algebraic curve with equation
\begin{equation*}
P(x) = (x_1^2+x_2^2)^2 + 18(x_1^2+ x_2^2) - 8x_1^3 + 24 x_1x_2^2 -27 = 0.
\end{equation*}

For this particular choice,  the matrix $(g_{ij})(x)$ is unique up to some scaling factor, and is  given by

\beq\label{matrix.real}\begin{cases}
g_{11}(x_1,x_2) = - (3x_1^2 - x_2^2 - 6x_1 - 9)\\
g_{12}(x_1,x_2) = -2x_2(2x_1+3)\\
g_{22}(x_1,x_2) = - (3x_2^2 - x_1^2 + 6x_1 - 9)
\end{cases} \eeq
Whence we deduce that $\det(g)=-3 P(x)$.
Moreover, in this representation, for the measure $\mu(dx)= c(\alpha)|P(x)|^\alpha dx$, the drift terms in the equation read
\begin{equation}
b_1(x_1,x_2) = -2( 6 \alpha+5)x_1, \quad b_2(x_1,x_2) = -2( 6 \alpha+5)x_2.
\end{equation}

The general operator $\LL^{(\alpha)}$ on the interior of Deltoid curve for which a family of orthogonal polynomial is formed of eigenvectors of $\LL^{(\alpha)}$  is therefore given by
\beqnas \LL^{(\alpha) }=&& g_{11}(x_1,x_2)\partial^2_{1}+ g_{22}(x_1,x_2)\partial^2_{2}+ 2  g_{12}(x_1,x_2) \partial^2_{1,2}\\&&  -2( 6 \alpha+5)x_1  \partial_{1}   -2( 6 \alpha+5) x_2\partial_ {2}\eeqnas      with associated measure $c(\alpha) \rho^\alpha dx$, with $$\rho=\frac{1}{3} \det(g)=  -(x_1^2+x_2^2)^2 -18(x_1^2+ x_2^2) +8 x_1^3 -24 x_1x_2^2 +27 .$$

 As long as we only deal with polynomials, it turns out that it is simpler to use complex variables. Indeed, let $Z= x_1+ix_2$ and  conjugate $\overline{  Z}= x_1-ix_2$, then  the generator is entirely characterized by
 \beq\label{matrix.complex}
  \begin{cases}\Gamma(Z,Z)=-4Z^2+12\bZ,\\
\Gamma(\bZ, Z)= -2Z\bZ+18,\\
\Gamma(\bZ, \bZ)= -4\bZ^2+12Z,\\
 \LL ^{(\alpha)}Z=  -2( 6 \alpha+5)Z , \LL^{(\alpha)} \bZ= -2( 6 \alpha+5)\bZ.\end{cases}
 \eeq

We can simplify the operator by  setting $Z= 3Z_1, \bar Z= 3\bar Z_1$ and multiplying $\LL^{(\alpha)}$ by $1/4$, which does not change the eigenvectors and multiply the eigenvalues by $1/4$. This  gives
\beq\label{matrix.complex2}
  \begin{cases}\Gamma(Z,Z)=\bZ- Z^2,\\
\Gamma(\bZ, Z)= 1/2(1-Z\bZ),\\
\Gamma(\bZ, \bZ)= Z- \bZ^2,\\
 \LL ^{(\alpha)}Z=  -1/2( 6 \alpha+5)Z , \LL^{(\alpha)} \bZ= -1/2( 6 \alpha+5)\bZ.\end{cases}
 \eeq

In particular,  giving a particular role to the case $\alpha= -1/2$, one has
\beq\label{comp.Lalpha} L^{(\alpha)}= L^{(-1/2)} -\frac{3}{2}(2\alpha+1)(Z\partial_Z + \bZ\partial_{\bZ}).\eeq
This model has been studied by \cite{Koorn1, Koorn2}, where the relationship with homogeneous spaces of rank 2 and root system $A_2$ has been put forward. Observe that the case $\lambda= -1/2$ corresponds to the Laplace-Beltrami operator associated with the Riemannian metric $g^{-1}$ associated with the inverse matrix of the matrix $g$.

Our aim here is to study the associated orthogonal polynomials together with the associated eigenvalues,  and various representations for it. Indeed, this family belongs to the larger class of Hall polynomials associated with root systems (here the root system $A_2$) (see \cite{macdonald}), and our aim here is to present some properties of these polynomials specific for this model.

\section{ $L^{(-1/2)}$ as a projection of the Euclidean Laplacian.\label{sec.laplacien2d}}

As already mentioned, the case $\alpha= -1/2$ corresponds to the Laplace-Beltrami operator associated to the inverse matrix $g^{-1}$.  If one computes the associated curvature (here, in dimension $2$, the scalar curvature is sufficient to characterize the metric), we may observe that it vanishes, and therefore it is not much surprising that the operator is the image, in the sense described in Section~\ref{sec1}, of the ordinary Laplace operator in $\R^2$.

In order to perform this identification, we first start by some remarks on the Deltoid curve.
\begin{figure}[ht]
 \centering		\includegraphics[width=.5\linewidth]{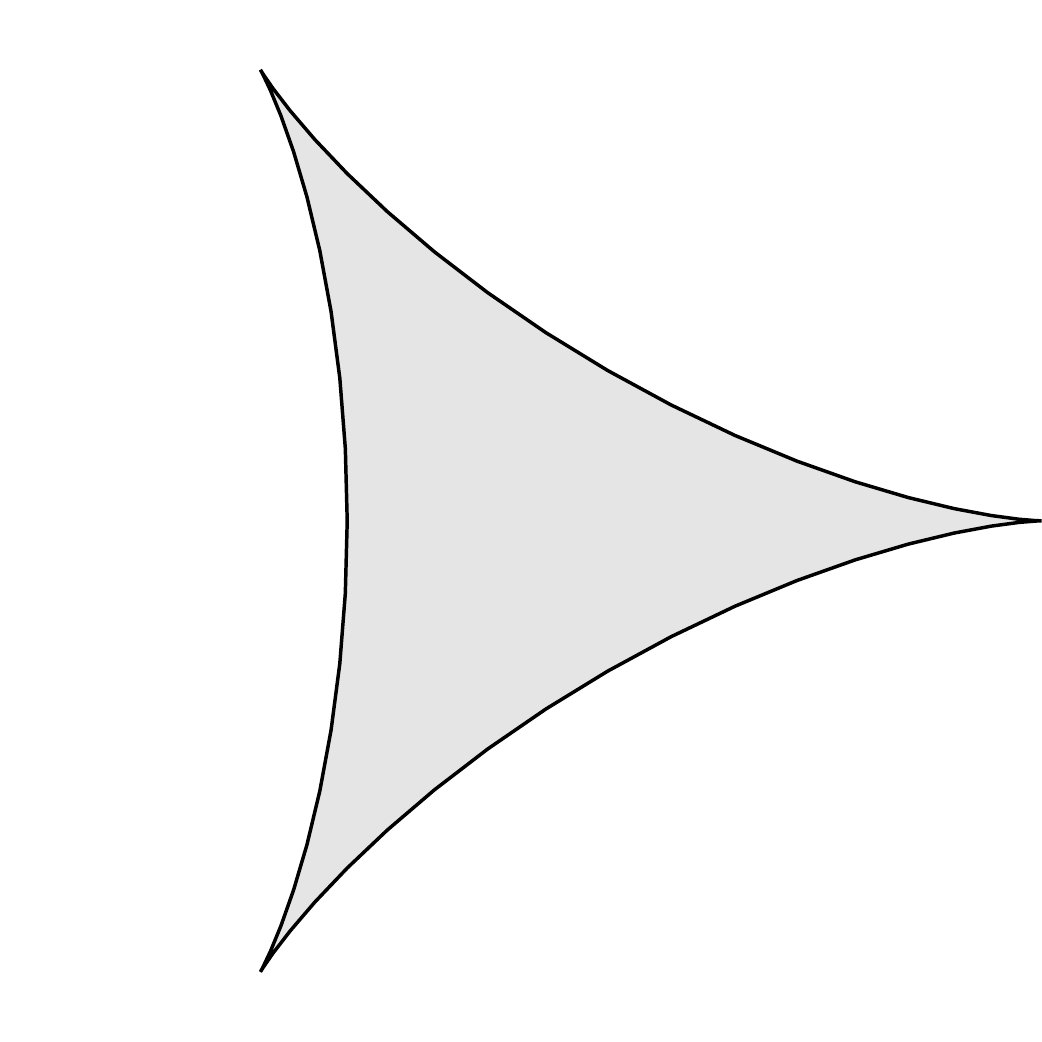}
		\caption{The deltoid domain.}
		\label{fig:Deltoide}
\end{figure}

It represented by the following parametric equations :
\begin{equation*}
x_1(\theta) = 2\cos\theta + \cos 2\theta, \quad x_2(\theta) = 2\sin \theta - \sin 2\theta
\end{equation*}
 In complex notations $z(\theta) = x_1(\theta) + ix_2(\theta) = 2e^{i\theta} + e^{-2i\theta}$.

We shall denote $\overline{\Omega}$ the interior of the deltoid curve, and our aim is to identify the above operator $\LL^{(-1/2)}$ as the image of the laplace operator $\Delta$ on $\R^2$ through the action of some function $Z : \R^2\mapsto \overline{\Omega}$ is the sense described above.

 To proceed, we shall first use the change of variables formula to see that
\begin{equation*}
\mathscr{\LL}^{(-1/2)}(f)(X,Y) = \Delta[f(X,Y)],
\end{equation*}
where $(X, Y):\R^2\mapsto \overline{\Omega} $, or in complex notations
\begin{equation*}
\mathscr{\LL}^{(-1/2)}(f)(Z,\overline{Z}) = \Delta[f(Z,\overline{Z})]
\end{equation*}
where
\begin{equation*}
Z = X + iY, \quad \overline{Z} = X - iY.
\end{equation*}

So that we are looking for some functions $Z : \R^2\mapsto \overline{\Omega}$ satisfying :
$$\Gamma_{\Delta}(Z,Z)=\nabla Z \cdot \nabla Z =-Z^2+\bZ,
\Gamma_{\Delta}(\bZ, Z)= \nabla \bZ \cdot \nabla Z =1/2(1-Z\bZ)$$
$$\Gamma_{\Delta}(\bZ, \bZ)=\nabla \bZ \cdot \nabla \bZ = -\bZ^2+Z,$$
$$ \Delta Z=  -Z , \Delta \bZ= -\bZ.$$

In what follows, for any $\beta\in \mathbb{R}^2\approx \mathbb{C}$,  $e_\beta(x)$ denotes the function  $\mathbb{R}^2\mapsto \mathbb{C}$, $x \mapsto e^{i(\beta\cdot x)}$

\bprop \label{Pro1}
Let $\beta_1, \beta_2, \beta_3 \in \mathbb{R}^2$ three vectors satisfying
\begin{equation*}
\beta_k\cdot \beta_k = 1, \quad \beta_k\cdot \beta_l = -\frac{1}{2},  l\neq k.
\end{equation*}
and let
\begin{equation*}
Z(x) = \sum_{k=1}^3 e_{\beta_k },
\end{equation*}
Then, one has
$$\Delta(Z)= -Z, ~\Delta(\bZ)= -\bZ, $$
and
$$\Gamma(Z, Z)= -Z^2+3\bZ , \Gamma(Z, \bZ)= \frac{1}{2}(9-Z\bZ), \Gamma(\bZ, \bZ)= -\bZ^2+3Z.$$
\eprop

\bpf Splitting $e_{\beta_k}(x), k \in \{1,2,3\}$ into a real and imaginary parts, one derives
\begin{equation*}
\left[\nabla e_{\beta_k} \cdot \nabla e_{\beta_l}\right]  = -(\beta_k \cdot \beta_l) e_{\beta_k+\beta_l}
\end{equation*}
for all $l,k \in \{1,2, 3\}$. As a result
\begin{align*}
\Gamma(Z,Z) &= -\sum_{k,l}(\beta_k\cdot\beta_l) e_{\beta_k+\beta_l} \\
& = - \left[\sum_{k=1}^3e_{2\beta_k} -\sum_{k < l} e_{\beta_k+\beta_l} \right]
 = - \left[Z^2 -3\sum_{k < l} e_{\beta_k+\beta_l} \right].
\end{align*}
But one easily checks that
\begin{equation*}
\sum_{k=1}^3\beta_k \cdot \sum_{k=1}^3\beta_k = 0
\end{equation*}
so that $\beta_1+\beta_2 + \beta_3 = 0$ yielding $\Gamma(Z,Z) = -Z^2 +3\overline{Z}$. Finally
\begin{align*}
\Gamma(Z,\overline{Z})  &=3- \sum_{k\neq l}^3 e_{\beta_k-\beta_l}
= 3-(\frac{Z\bZ-3}{2}).
\end{align*}
The identification is obtained changing $(Z, \bZ)$ into $(Z/3, \bZ/3)$.
\epf

\begin{Nota}
The assumption $\beta_k \cdot \beta_l = -1/2, k \neq l$ is by no means a loss of generality. Indeed, if one rather assumes that $\beta_k \cdot \beta_l = c \in (-1,1)$ then one may take $\beta_1 = 1$ due to the rotation invariance of our conditions. But then $\beta_1 \cdot \beta_2 = \beta_1 \cdot \beta_3$ forces both first coordinates of $\beta_2,\beta_3$ to be equal $c$, while the fact that the vectors have unit length entails
\begin{equation*}
\beta_2 = c + i\sqrt{1-c^2}, \quad  \beta_3 = c - i \sqrt{1-c^2}.
\end{equation*}
Together with $\beta_2 \cdot \beta_3 = c$ show that $c$ to be a root of $2c^2 - c - 1$ yielding finally $c = -1/2$. Hence $\beta_1 = 1, \beta_2 = j, \beta_3 = j^2$, the cubic roots of the unit, up to an orthogonal transformation.
\end{Nota}

From now on, with no loss of generality, we shall assume that  $(\beta_1, \beta_2, \beta_3) = (1,j, j^2)$.

 One immediately sees that $Z$ is invariant under the action (by translation) of the lattice $\mathbb{L}$ generated by $4\pi \beta_1, 4 \pi \beta_2$,
\begin{equation*}
\mathbb{L}=  4\pi \mathbb{Z}\beta_1 + 4\pi \mathbb{Z}\beta_2,
\end{equation*}
by rotation with $2\pi/3$ angles and by symmetry with respect to the lines $\mathbb{R}\beta_k$.  One may also observe that it is invariant under the symmetry with respect to the horizontal line $\{y= 2\pi/\sqrt{3}\}$.
From this, one  sees that $Z$ is also  invariant under the symmetries with respect of the lines of the regular  triangular lattice  $\bbL_1$ whose fundamental domain is the regular triangle $\mathcal{A}$ whose vertices are $(0,0),(4\pi/3,0), (4\pi/3)e^{i\pi/3}$ (see below). We shall say in the sequel that a function having those invariance  have the symmetries of the lattice $\bbL_1$.

Equivalently, $Z$ is invariant under the action of the dihedral group $\mathcal{D}_3$ of affine type (\cite{hum}).  As a matter of fact, $Z$ is uniquely determined by its restriction to  $\mathcal{A}$.

\medskip
\begin{picture}(100,180)(1,1)
\begin{thinlines}
\multiput(10,30)(30,0){9}{\line(3,5){90}}
\multiput(130,30)(30,0){9}{\line(-3,5){90}}
\multiput(10,30)(0,25){7}{\line(1,0){360}}
\put(279,30){\line(3,5){90}}
\put(310,30){\line(3,5){60}}
\put(340,30){\line(3,5){30}}
\put(369,80){\line(-3,5){60}}
\put(369,130){\line(-3,5){30}}
\put(100,30){\line(-3,5){90}}
\put(70,30){\line(-3,5){60}}
\put(40,30){\line(-3,5){30}}
\put(10,80){\line(3,5){60}}
\put(10,130){\line(3,5){30}}
\put(130,0){\textrm{\emph{The lattice of regular triangles $\bbL_1$.}}}
\put(185,110){$\mathcal{A}$}
\put(172,92){$0$}
\put(198,92){${\tiny4\pi/3}$}
\end{thinlines}
\end{picture}
\medskip

\bprop
$Z$ is a one-to-one map from $\partial \mathcal{A}$ onto $\mathcal{D} = \partial \Omega$ and from $\overline{\mathcal{A}}$ onto $\overline{\Omega}$.  In particular, it maps the whole plane onto $\overline{\Omega}$.
\eprop
\bpf Recall the parametric equation of $\mathcal{D}$ written in complex notations
\begin{equation*}
z(\theta) = 2e^{i\theta} + e^{-2i\theta}
\end{equation*}
so that $z(\theta) = Z(-2\theta,0)$, where $\theta$ runs over any interval of length $2\pi$. Then the invariance of $Z$ under rotations of angles $\pm 2\pi/3$ shows that the images of the intervals
\begin{equation*}
[-4\pi/3,0], \quad [4\pi/3, 8\pi/3]
\end{equation*}
coincides with the images of the oblique edges of $\mathcal{A}$ (the cusps of $\mathcal{D}$ are the images of $\{\theta= 0, 4\pi/3, 8\pi/3$\}). Thus $Z$ maps $\partial{A}$ onto $\mathcal{D}$ and it is easy to check from the complex parametrization of $\mathcal{D}$ that $Z$ is one-to-one there. Combined with the compactness of $\mathcal{A}$ and the continuity of $Z$, we deduce that $Z$ maps $\mathcal{A}$ into $\overline{\Omega}$.
But then $Z(\mathcal{A}) = \overline{\Omega}$ since otherwise $\overline{\Omega}$ would not be simply connected, which leads to a contradiction.
Now, we shall use the following parametrization of $\overline{\Omega}$:
\begin{equation*}
Z(x_1,x_2)= e^{ix_1} + 2e^{-i\frac{x_1}{2}}\cos(\frac{\sqrt{3}}{2}x_2), \, x=(x_1,x_2) \in \mathcal{A}.
\end{equation*}
For fixed $x_1 \in [0, 2\pi/3]$, the image by $Z$ of the vertical segments
\begin{equation*}
[(x_1,0), (x_1, \sqrt{3}x_1)] \in \mathcal{A}
\end{equation*}
is the line segment $I(x_1) = [A(x_1),B(x_1)]$ where
\begin{eqnarray*}
A(x_1) &=& (\cos(x_1)+2\cos(\frac{x_1}{2}), \sin(x_1) - 2\sin (\frac{x_1}{2})) \\
B(x_1) &=& (2\cos(x_1)+\cos(2x_1), 2\sin(x_1) - \sin(2x_1)).
\end{eqnarray*}
Thus, the coordinates of $A(x_1)$ are decreasing as functions of the variable $x_1$ while those of $B(x_1)$ are decreasing and increasing respectively. Equivalently, $A(x_1)$ runs over the half part of the lowest branch of $\mathcal{D}$ starting from $(3,0)$ while $B(x_1)$ runs over the whole highest one since clearly $B(x_1) = A(-2x_1)$. As a matter of fact, two line segments $I(x_1), I(x_1'), 0 \leq x_1 \neq x_1' \leq 2\pi/3$ never intersect. A similar reasoning applies when $x_1 \in [2\pi/3,4\pi/3]$ and the segment $[x_1, -\sqrt{3}x_1 + 4\pi/\sqrt{3}]$, $A(x_1)$ runs over the remaining half part of the lowest branch while $B(x_1)$ runs over the whole third one. As a matter of fact, $Z$ is a one-to-one from $\overline{\mathcal{A}}$ onto $\overline{\Omega}$. Finally, $\overline{\mathcal{A}}$ is a fundamental domain for the action of the affine group $\mathcal{D}_3$ on $\mathbb{R}^2$ so that every $x \in \mathbb{R}^2$ is conjugated to a unique element of $\mathcal{A}$. The proposition is proved. $\hfill \blacksquare$

 \begin{figure}[ht]

 \centering		\includegraphics[width=.3\linewidth]{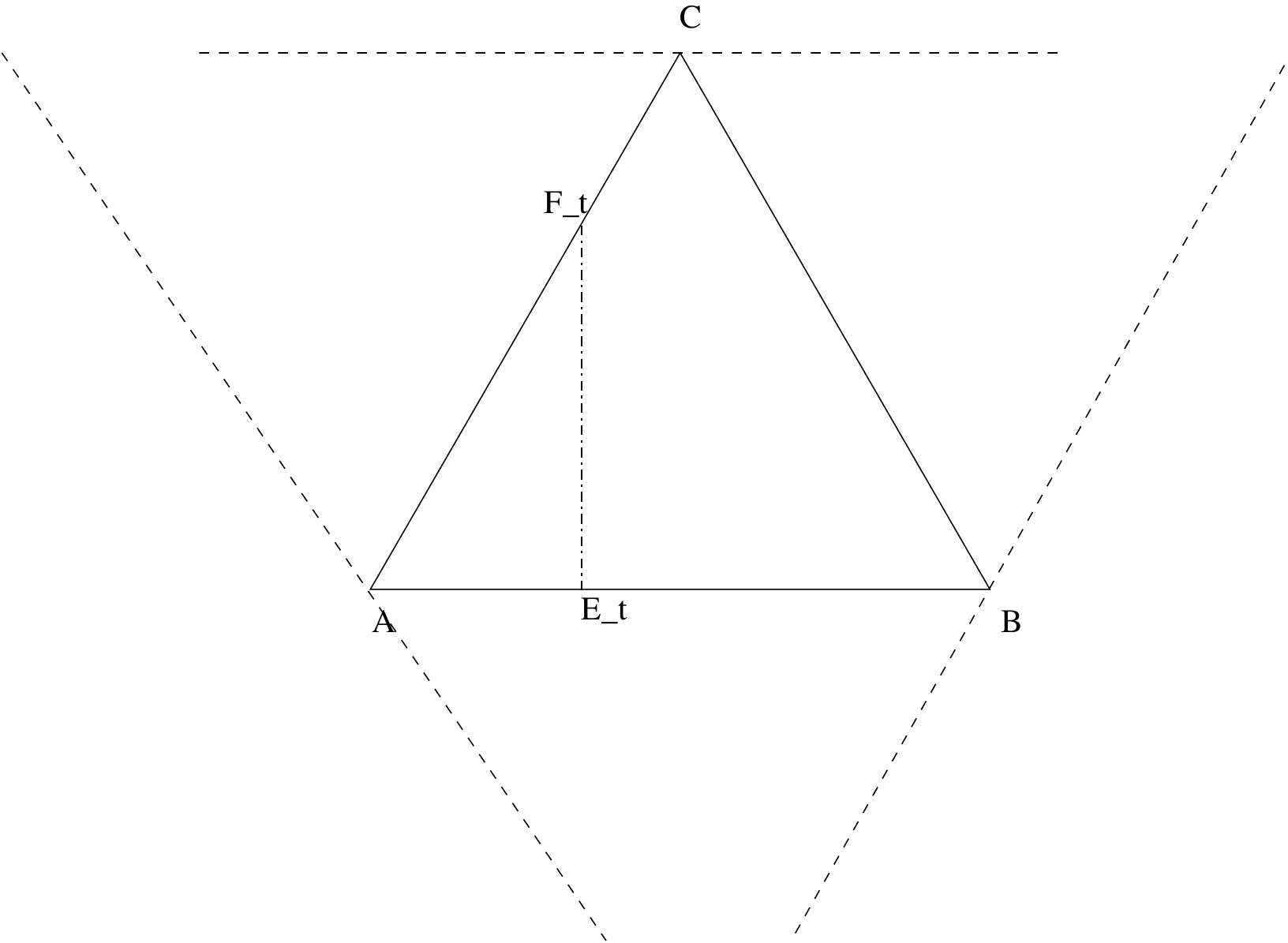}

		\label{fig:triangleABC}
\end{figure}
\epf

We are now in situation to identify  the operator $\LL^{(-1/2)}$  as an image of the 2-dimensional Laplace operator acting on functions which are invariant under the symmetries of the lines in the triangular  lattice.

\bprop A measurable function $f : \bbR^2\mapsto \bbR$  have the symmetries of the lattice $\bbL_1$ if and only if it may be written $f= g(Z)$, where $g: \Omega \mapsto \bbR$ is a measurable function.

Moreover, when $f\in \cC^2$, then we may chose $g\in \cC^2$, in which case
\beq\label{change.var.Z}\Delta(g(Z))= \LL^{(-1/2)} (g)(Z).\eeq
In other words, $\LL^{(-1/2)}$ is nothing else than the 2-dimensional Laplace operator acting on set of functions having the symmetries of the lattice $\bbL_1$.

\eprop

\bpf If we denote by $Z^{-1}$ the inverse map $\Omega \mapsto A$ of the restriction of $Z$ to $A$, then we just set $g= f\circ Z^{-1}$.

Moreover, the change of variable formula~\eqref{change.var}  and formulae given in Proposition~\ref{Pro1} give immediately~\eqref{change.var.Z}.

\epf

\brmq\label{rmqdetvdm}

It is worth to remark that if we set  that  $Z=z_1+z_2+ z_3$, where $z_i$ are complex numbers such that $|z_i|=1$ and $z_1z_2z_3=1$, then  $$ \rho=-3(z_1-z_2)^2(z_2-z_1)^2(z_3-z_1)^2  .$$   \\
 Indeed, if $\rho(x,y)$ is the determinant of the matrix $(g^{ij})$ written in $(x_1,x_2)$ coordinates, one sees that, in $(Z, \bar Z)$ coordinates, it may be written as
 $$\frac{1}{4}\Big(\Gamma(Z, \bar Z)^2-\Gamma(Z,Z)\Gamma(\bar Z, \bar Z)\Big),$$ which gives
\beq
\rho=12(Z^3+\bZ^3)-3Z^2\bZ^2-54Z\bZ+81
\eeq
\ermq

Using remark \ref{rmqdetvdm} and the above diffeomorphism between the deltoid and the triangle, one obtains
\bprop\label{prop.int.mesure} The function  $\det(g)^\alpha $ on the deltoid is integrable with respect to the Lebesgue measure if and only if $\alpha> -5/6$.

\eprop
In the sequel, we shall always set $\lambda=\frac{1}{2} (6\alpha+5)$.

\bpf
By the change variables formula we have
\beqnas
\int_{\mathcal{D}} \det(g)^{\alpha} dx_1dx_2=\int_{ \mathcal{A}} \det(g)^{\alpha+\frac{1}{2}} dx_1dx_2\\
= (-3)^{\alpha+\frac{1}{2}}\int_{ \mathcal{A}}( (z_1-z_2)(z_2-z_3)(z_3-z_1))^{2\alpha+1} dx_1dx_2
\eeqnas
where $z_1=e^{i\theta_1}, z_2=e^{i\theta_2}, z_3=e^{i\theta3}$  and $\theta_1=x_1,\theta_2=-\frac{x_1}{2}+\frac{\sqrt3 x_2}{2}, \theta_3=-(\theta_1+\theta_2)= -\frac{x_1}{2}-\frac{\sqrt3 x_2}{2} $\\

A rapid inspection of the  integrability condition for this function on the triangle shows that, near the boundary and  outside the corners of the triangle, the integrability  condition is $\alpha>-1$, while at the corner of the triangles, the condition is more restrictive.
Indeed, for the integrability of the measure near the  point $(0,0)$,  then
$(z_1-z_2)(z_2-z_3)(z_3-z_1)) \simeq -i\frac{3}{4}\sqrt{3}x_2(x_2^2-3x_1^2)$ and if we set $x_2=\sqrt3tx_1, t \in [0,1]$ we have
$\det(g)^{\alpha +\frac{1}{2}}\simeq \big(-\frac{27}{4}t^2x_1^6(1-t^2)^2\big)^{\alpha +\frac{1}{2}}$, which is  integrable for the measure $tdtdx_1$ if and only if $\alpha > -\frac{5}{6}$.

\epf

\section{The $\LL^{(1/2)}$ as a projection of the Casimir operator on $SU(3)$\label{sec.su3}}~

Let $\cG$ be a compact semi simple Lie linear group with Lie algebra $\cL$ , seen at the tangent space at Id for $G$, with Lie bracket $[A,B]$ (see \cite{Faraut}). On $\cL$, the Killing form is  a scalar product defined by $\lag A, B\rag= -\tr (A.B)$
. On the other hand, to any $A\in \cL$ is associated a vector field $X_A$ on $\cG$ defined as
$X_A(f)(g)= \partial_t\mid_{t=0} f(ge^{tA})$. Given an orthonormal basis $(A_1, \cdots, A_d)$ in $\cL$ with respect to the Killing form, the  Casimir operator is defined as $\Delta_G=\sum_i X_{A_i}^2$.  It is also the Laplace-Beltrami operator on $\cG$ when $\cG$ inherits the Riemannian structure from the Killing form in $\cL$. $\Delta_\cG$ is a second order differential operator in the sense that it satisfies the change of variable formula~\eqref{change.var}. We shall denote by $\Gamma_\cG$ the corresponding carr\'e du champ operator.

The Casimir operator commutes with the Lie group action. More precisely, if, for $g\in \cG$ and  for any function $f : \cG\mapsto \R$,  one defines the right action $R_g(f)(k)= f(kg)$, then $\Delta_\cG R_g= R_g \Delta_\cG$, and the same holds true for the left action $L_g(f)(k)= f(gk)$.

In order to entirely determine the action of $\Delta_\cG$ on functions of $\cG$, it is enough to compute $\Delta_\cG (f_i)$ and $\Gamma_\cG(f_i,f_j)$ for a set of functions which generates all functions on $\cG$ (say as $\sigma$-algebras). Once again, it could be helpful to consider complex valued functions, and  on $SU(n)$, if one represents $g$ as a matrix $(z_{ij})$ with complex entries, we shall consider the coordinates $g\mapsto z_{ij}$ and $g\mapsto \bz_{ij}$ as generating functions.

When performing the above computations in $SU(n)$, one ends up with the following formulae
\beqnas&& \Delta_{SU(n)}(z_{kl})= -2\frac{(n-1)(n+1)}{n} z_{kl},~\Delta_{SU(n)}(\bar z_{kl})= -2\frac{(n-1)(n+1)}{n}\bar z_{kl},\\
&&\Gamma_{SU(n)}(z_{ij}, z_{kl})= -2z_{il}z_{kj}+\frac{2}{n}z_{ij}z_{kl},~\Gamma_{SU(n)}(\bz_{ij}, \bz_{kl})= -2\bz_{il}\bz_{kj}+\frac{2}{n}\bz_{ij}\bz_{kl}, \\
&&\Gamma_{SU(n)}(z_{ij}, \bz_{kl})= 2(\delta_{ik}\delta_{jl} -\frac{1}{n}z_{ij}\bz_{kl}).
\eeqnas

For any $p\in \Z$, consider  the functions $SU(n)\mapsto \C$: $T_p(g)= \tr(g^p)$.  For $p\geq 1$, one has
\beq\label{expand.trace}T_p(g) = \sum_{i_1, \cdots, i_p= 1}^n z_{i_1i_2}z_{i_2i_3} \cdots z_{i_pi_1},\eeq while the same formula holds for $p\leq -1$ replacing $z_{ij}$ by $\bz_{ij}$ (and of course $T_{-p}= \overline{T}_p$).
From the change of variable formula, one has, for any $m$-uple of functions $(f_1, \cdots, f_m)$ and any diffusion generator $\LL$
\beq\label{formule.produit}\LL(f_1\cdots f_m)=  f_1\cdots f_m\Big(\sum_{i=1}^m \frac{\LL f_i}{f_i} + \sum_{i,j=1}^m\frac{ \Gamma(f_i,f_j)}{f_if_j} - \sum_{i=1}^m\frac{ \Gamma(f_i,f_i)}{f_i^2}\Big),\eeq
and, for any $m$-ulple $(f_1, \cdots, f_m)$ and any $k$-uple $(g_1, \cdots, g_k)$
$$\Gamma(f_1\cdots f_m, g_1\cdots g_k)= f_1\cdots f_m g_1\cdots g_k\Big(\sum_{i=1}^m\sum_{j=1}^k \frac{\Gamma(f_i, g_j)}{f_ig_j}\Big).$$

Applying these to the explicit expression \eqref{expand.trace} of  $T_p$, one gets, for $p\geq 1$
\beq\label{eqLTpSU3}\Delta_{SU(n)} T_p = -p\Big(2 (\frac{n^2-p}{n})T_p + \sum_{i=1}^{p-1} T_iT_{p-i}\Big),\eeq with the conjugate formula for $p\leq -1$,  while, for any $p,q\in \Z$
\beq\label{eqGammaTpSU3}\Gamma_{SU(n)}(T_p, T_q) = 2|pq|\big( \frac{T_pT_q}{n}-T_{p+q}\big).\eeq

In particular, if we set $Z= T_1$, $\bZ= T_{-1}$. Then
\beq\label{delta.sun}\Delta_{SU(n)} Z = -2 \frac{n^2-1}{n}Z, \Delta_{SU(n)} \bZ=-2 \frac{n^2-1}{n}\bZ,\eeq
and
\beqna \label{gamma.sun}&\Gamma_{SU(n)}(Z,Z)= 2(\frac{Z^2}{n}-T_2), ~\Gamma_{SU(n)}(\bZ,\bZ)= 2(\frac{\bZ^2}{n}-\overline{T}_2), \\
&\Gamma_{SU(n)}(Z,\bZ))= 2(3-\frac{Z\bZ}{n}).\eeqna

Now, consider more precisely the case $n=3$. For any matrix in $SU(3)$, if $(\mu_1, \mu_2, \mu_3)$ denote its eigenvalues, $T_p= \mu_1^p+\mu_2^p+ \mu_3^p$.  The $\mu_i$ are complex numbers with $|\mu_i|=1$ and $\mu_1\mu_2\mu_3=1$. With $Z= \mu_1+\mu_2+\mu_3$, they are solution of the equation
$$X^3-ZX^2+ \bZ X-1=0,$$ and multiplying this by  $X^p$ and summing over the three values $\mu_1, \mu_2, \mu_3$, one gets for any $p\in \Z$
\beq\label{recurrence.trace.su3}T_{p+3} -ZT_{p+2} +\bZ T_{p+1}-T_p= 0,\eeq which, for $p=-1$ gives
$T_2= Z^2-2\bZ$, and similarly $\overline{T}_2= \bZ^2-2Z$. Replacing this values in~\eqref{delta.sun} and~\eqref{gamma.sun} lead to
$$\Delta_{SU(3)} Z= -\frac{16}{3}Z, ~\Delta_{SU(3)}(\bZ)=  -\frac{16}{3}\bZ$$ and
\beqnas &&\Gamma_{SU(3)}(Z,Z)= \frac{4}{3}(3\bZ-Z^2), \Gamma_{SU(3)}(\bZ,\bZ)= \frac{4}{3}(3Z-\bZ^2), \\&&\Gamma_{SU(3)}(Z,\bZ)= \frac{2}{3}(Z\bZ-9).\eeqnas
It remains to replace $Z$ by $Z/3$ to observe that $\frac{3}{4}\Delta_{SU(3)}$, acting on functions of $(Z, \bZ)$ is nothing else than $\LL^{(1/2)}$ and  $\lambda=4$. Observe also that the functions on $SU(3)$ which depend only on $(Z, \bZ)$ are exactly those functions which depend only on the spectrum of the matrix $g\in SU(3)$, that is the functions which are invariant under $g\mapsto h^{-1}gh$, for any $h\in SU(3)$. Indeed, as long as polynomials are concerned,  those functions are exactly functions depending only on the traces $T_p, p\in \Z$, and formula~\eqref{recurrence.trace.su3} shows that these functions are again polynomials in the variables $(Z, \bZ)$.

\section{Eigenvalues and eigenvectors\label{sec.eigen} }
We proceed now to the determination of the  eigenvalues of $\LL^{(\alpha)}$, and give a recurrence formula for the corresponding eigenvectors.

In dimension 1, it is well known (and easy to check) that for any probability measure for which the polynomials are dense in $\cL^2(\mu)$, the unique (up to the sign) associated sequence of orthogonal polynomials satisfies a 3 term recurrence formula (see \cite{NSU}), usually written under the form

$$xP_n = a_nP_{n+1} + b_n P_n + a_{n-1} P_{n-1}.$$
It is not the case in dimension $2$, since one then would get in general a recurrence formula involving at each step $n$ an increasing number of terms.  Indeed, if, for each degree $n$,  one denotes by $\cP_n$ the space of polynomials of total degree less than $n$ and by $\cH_n$  the space of polynomials  in $\cP_n$ orthogonal to $\cP_{n-1}$, one has for any polynomial orthogonal polynomial $P\in \cH_n$
$$x_1P= Q_{n+1} + Q_n + Q_{n-1}, ~x_2 P= R_{n+1} + R_n + R_{n-1},$$ where $Q_i$ and $R_i$ belong to $\cH_i$. But the spaces $\cH_i$ have dimension $i+1$, and one should in general not  expect any simple  recurrence formula.

However, looking more precisely at the form of the operators $\LL^{(\alpha)}$ is the variables $(Z, \bZ)$, one should  expect for the sequence of eigenvectors of $\LL^{(\alpha)}$ a 6 term recurrence formula. It comes as a surprise that indeed one is able to get a 3 term recurrence formula as in dimension 1.

 We first start by investigating the eigenvalues.   Recall first that we are looking for polynomials $P_{p,q}^{(\alpha)}$ such that
\begin{equation*}
\LL^{(\alpha)} (P_{p,q}^{(\alpha)}) = -\lambda_{p,q} P_{p,q}^{(\alpha)}
\end{equation*}
where $ (p,q) \in \mathbb{N}^2$ is a bi-index whose weight $p+q$ is the degree of $P^{(\alpha)}_{p,q}$.

\bprop The eigenvalues of $\LL^{(\alpha)}$ are
$$\lambda_{p,q}= (\lambda-1)(p+q)+ p^2+q^2+ pq,$$where $\lambda= \frac{1}{2}(6\alpha+5)$.

\eprop

\bpf
The complex representation easily leads to the eigenvalues. Indeed, if $\cP_n$ denotes the space of polynomials (now in  the variables $(Z, \bZ)$) with total degree at most $n$, one may write any $P\in \cP_n$ as
$$P= \sum_{p=0}^n a_{p,q} Z^p\bZ^{n-p} + Q= P_n + Q,$$ where $Q\in \cP_{n-1}$.

Now, looking at the action of $\LL^{(\alpha)}$ on the highest degree term $P_n$  of $P$, one sees that if $\LL^{(\alpha)} P= -\mu P$, then the highest degree term  $\hat P_n$ of $\LL^{(\alpha)} P_n$ is equal to $-\mu \hat P_n$. It remains to observe the action of $\LL^{(\alpha)}$ on those highest terms. Fortunately, in coordinates $(Z, \bZ)$, this action is diagonal (which is not the case in coordinates $(x_1,x_2)$).

Indeed, the change of variable formula \eqref{change.var}   gives
\beqnas \LL^{(\alpha)} (Z^p\bZ^q)=&& pZ^{p-1}\bZ^q\LL^{(\alpha)} Z
\\&&
+ qZ^{q-1}Z^p\LL^{(\alpha)} \bZ + p(p-1)Z^{p-2}Z^q\Gamma(Z,Z)
\\&&
+ 2pqZ^{p-1}\bZ^{q-1}\Gamma(Z,\bZ)+ q(q-1)\bZ^{q-2}Z^p\Gamma(\bZ,\bZ),\eeqnas whose highest term is
$$-\lambda_{p,q} Z^p\bZ^q,$$
with $\lambda_{p,q}= (\lambda-1)(p+q)+ p^2+q^2+ pq$, where $\lambda= \frac{1}{2}(6\alpha+5)$.

\epf
\brmq
When $\alpha\notin \Q$, then the eigenspaces associated to the eigenvalues $\lambda_{p,q}$ are at most two-dimensional (and exactly two dimensional when $p\neq q$). Indeed, writing $\sigma= p+q$ and $\pi= pq$, with similar notation $\pi',\sigma'$ for $(p',q')$, we have
$\lambda_{p,q}= (\lambda-1)\sigma+\sigma^2-\pi$, and therefore if $\lambda_{p,q}= \lambda_{p',q'}$, then either $\sigma=\sigma'$ and then $\pi= \pi'$, either
$\dis \lambda= 1-\sigma-\sigma'+ \frac{\pi-\pi'}{\sigma-\sigma'}$, whence $\lambda\in \Q$.
\ermq

We shall see moreover that for any $(p,q)$ there exists exactly one polynomial $P^{(\alpha)}_{p,q}(Z, \bZ)= Z^p\bZ^q+ \hbox{\rm lower degree terms}$ which is an eigenvector of $\LL^{(\alpha)}$.
Indeed,
\bthm\label{defPpq} Define the family of polynomials   $P^{(\alpha)}_{p,q}(Z, \bZ)$ by induction from
$$P^{(\alpha)}_{0,0}= 1~, P^{(\alpha)}_{0,1}= \bZ~, P^{(\alpha)}_{1,0}=Z, $$
and
\beq\label{defPpq2}\begin{cases}P^{(\alpha)}_{p+1, q}&= ZP^{(\alpha)}_{p,q}+a_{1}(\lambda,p)P^{(\alpha)}_{p-1,q+1} +a_{2}(\lambda,p,q)P^{(\alpha)}_{p,q-1},\\
\\
P^{(\alpha)}_{p, q+1}&= \bZ P^{(\alpha)}_{p,q} +a_{1}(\lambda,q)P^{(\alpha)}_{p+1,q-1} +a_{2}(\lambda ,q,p)P^{(\alpha)}_{p-1,q} ,\end{cases}\eeq
where \\

\beqnas\begin{cases} a_{1}(\lambda,p)&= {\displaystyle -\frac{ p(3p+2\lambda-5)} {(\lambda+3p-1)(\lambda+3p-4)}}\\
\\
~ a_{2}(\lambda,p,q)&={\displaystyle -\frac{N_{p,q}}{D_{p,q}}}~ ~ ~ ~
\end{cases}
\eeqnas

where
\beqnas \begin{cases}
N_{p,q}&=q(3q+2\lambda-5)(\lambda+3(p+q)-1) (\lambda+p+q-2) \\
\\
~D_{p,q}&=(\lambda+3q-1)(2\lambda+3(p+q)-5) (2\lambda+3(p+q)-2) (\lambda+3q-4)
\end{cases}
\eeqnas

\beq
 \lambda=\frac{1}{2}(6\alpha+5)>0
 \eeq
 \brmq The only possible values of $\lambda$ for which the denominators vanishes in the above formulae are $\lambda= 1$ and $\lambda = 4$, which correspond to the $(p,q)\in \{(0,0),(1,0),(0,1)\}$. In those situations, we have to replace $a_1(\lambda,p)$ and $a_2(\lambda,p,q)$  by :
  \beqnas
   \begin{cases}
   a_1(\lambda, p)&= \lim_{\epsilon \to 0} a_1(\lambda +\epsilon, p)\\
   \\
   a_2(\lambda,p,q)&=\lim_{\epsilon \to 0} a_2(\lambda +\epsilon, p, q)
   \end{cases}
   \eeqnas

Moreover, $a_1(1,p)= a_1(4,p)= -1/3$  and $a_2(1,p,q)= a_2(4,p,q)= -1/9$ for every $(p,q)$ except for those values of $(p,q) \in \{(0,0),(1,0),(0,1)\}$ .
We have indeed

\beqnas
a_1(1,1)&=& -\frac{2}{3},a_2(1,0,1)=-\frac{1}{3} ,\\
a_1(4,1)&=& -\frac{1}{3},a_2(4,0,1)= -\frac{1}{9}, \\
a_1(\lambda,0)&=&a_2(\lambda,p,0)=0,\eeqnas

In the case $\alpha=1/2$ the recurrence formulae  simplify and for every
$p,q \ge 0 $ is :

\beq
P^{(1/2)}_{p+1, q}= ZP^{( 1/2)}_{p,q}-\frac{1}{3}P^{(1/2)}_{p-1,q+1} -\frac{1}{9}P^{( 1/2)}_{p,q-1},\\
\eeq

But in the other case $\alpha=-1/2$, the recurrence formulae is the same  except for the values  $(p,q)=\{(1,0), (0,1)\}$ corresponding to the polynomials $P^{(-1/2)}_{1,1}=Z\bZ-\frac{1}{3}$, and
 $P^{(-1/2)}_{2,0}=Z^2-\frac{2}{3}\bZ$.

Therefore, for $\alpha=-1/2, 1/2$, the recurrence formulae for  the polynomials are the same, except for the first two coefficients.

  \ermq

 Then, for  the operator $\LL^{(\alpha)}$ determined from~\eqref{matrix.complex2}, we have
$$L^{(\alpha)}P^{(\alpha)}_{p,q}=-\lambda_{p,q}P^{(\alpha)}_{p,q}$$ where $$\lambda_{p,q}=(\lambda-1)(p+q)+p^2+q^2+pq.$$
\ethm
It is worth to observe that since $a_1(\lambda,0)=a_2(\lambda,p,0)=0$, formula~\eqref{defPpq2} make sense for $p=0$ and $q=0$, and defines completely the family $P_{p,q}^{(\alpha)}$ for any $(p,q)\in \N^2$.
One observes that  $P^{(\alpha)}_{p,q}(Z, \bZ)= Z^p\bZ^q+ \hbox{\rm lower degree terms}$, and have real coefficients.  It is also easily checked that $P^{(\alpha)}_{q,p}= \overline{P^{(\alpha)}_{p,q}}$.

The proof of Theorem~\ref{defPpq} is rather tedious. We start with a Lemma:
\blem\label{lemmaGammaPpq} For the same family of polynomials defined in~\eqref{defPpq2} and the $\Gamma$ operator defined in~\eqref{matrix.complex2}, we have
$$\Gamma(Z,P^{(\alpha)}_{p,q})= \alpha_{0}(p,q)P^{(\alpha)}_{p+1,q}+\alpha_{1}(p,q)P^{(\alpha)}_{p-1,q+1} +\alpha_{2}(p,q)P^{(\alpha)}_{p,q-1}$$
$$\Gamma(\bZ,P^{(\alpha)}_{p,q})= \alpha_{0}(q,p)P^{(\alpha)}_{p,q+1}+\alpha_{1}(q,p)P^{(\alpha)}_{p+1,q-1} +\alpha_{2}(q,p)P^{(\alpha)}_{p-1,q}$$
where
 \beqnas
  \alpha_{0}(p,q)&=&{\displaystyle-\frac{1}{2}(q+2p)},\\
  \\
   \alpha_{1}(p,q)&=&{\displaystyle\frac{1}{2}\frac{p(2\lambda +3q-5)(\lambda+p-q-1)}{(\lambda+3p-1)(\lambda+3p-4)}} \\
   \\
\alpha_{2}(p,q)&=&\frac{1}{2} \frac{N1_{p,q}}{D1_{p,q}}
\eeqnas
where
$$\begin{cases}N1_{p,q}=q(3q+2\lambda-5)(\lambda+3(p+q)-1) (\lambda+p+q-2)(2\lambda+p+2q-2).\\
 D1_{p,q}= (\lambda+3q-1)(2\lambda+3(p+q)-5) (2\lambda+3(p+q)-2) (\lambda+3q-4).
 \end{cases}
 $$

\elem
It is worth to observe that although the definition of $\Gamma$ does not involve the parameter $\alpha$ (or equivalently the parameter $\lambda$), the recurrence formula defining $P^{(\alpha)}_{p,q}$ does, and this Lemma is valid whatever the parameter $\alpha$ is. However, it is not clear from the proof below for which family of  recurrence formulae on $P^{(\alpha)}_{p,q}$   three terms recurrence formulae for $\Gamma(Z, P^{(\alpha)}_{p,q})$ and $\Gamma(\bZ, P^{(\alpha)}_{p,q})$ are  still valid.
\begin{proof}
In what follows, we remove the parameter $\alpha$ from the formulae, since it shall not change up to end of this Section.
 Lemma~\ref{lemmaGammaPpq}   proved by induction, from
\beqnas
 \Gamma(Z,P_{p+1,q})&=&\Gamma(Z,ZP_{p,q})+a_{1}(\lambda, p)\Gamma(Z,P_{p-1,q+1})\\
 &&+a_{2}(\lambda,p,q)\Gamma(Z,P_{p,q-1})\\
  &=& \Gamma(Z)P_{p,q}+Z\Gamma(Z,P_{p,q})+ a_{1}(\lambda,p)\Gamma(Z,P_{p-1,q+1})\\&&
+ a_{2}(\lambda,p,q)  \Gamma(Z,P_{p,q-1})\eeqnas
 and finally
\beqnas\label{recGamma} \Gamma(Z,P_{p+1,q}) &=&( \bZ-Z^2)P_{p,q}+Z\Gamma(Z,P_{p,q})\\&&
+ a_{1}(\lambda,p)\Gamma(Z,P_{p-1,q+1})+ a_{2}(\lambda,p,q)  \Gamma(Z,P_{p,q-1})
\eeqnas
And  by the definition \eqref{defPpq} we have
\beqnas
\bZ P_{p,q}&= &P_{p, q+1} -a_{1}(\lambda,q)P_{p+1,q-1} -a_{2}(\lambda,q,p)P_{p-1,q}  \\
Z P_{p,q}&= &P_{p+1, q} -a_{1}(\lambda,p)P_{p-1,q+1} -a_{2}(\lambda,p,q)P_{p,q-1} \\
\eeqnas
So that
\beq
  Z^{2}P_{p,q}=Z(ZP_{p,q})= ZP_{p+1, q} -a_{1}(\lambda, p)ZP_{p-1,q+1}-a_{2}(\lambda, p,q)ZP_{p,q-1}\eeq

Furthermore,
  \beqnas\begin{cases}
 ZP_{p+1, q} =P_{p+2,q}-a_{1}(\lambda, p+1)P_{p,q+1}-a_{2}(\lambda,p+1,q)P_{p+1,q-1}\\
 ZP_{p-1,q+1}=P_{p,q+1}-a_{1}(\lambda, p-1)P_{p-2,q+2}-a_{2}(\lambda,p-1,q+1)P_{p-1,q}\\
 ZP_{p,q-1}= P_{p+1, q-1} -a_{1}(\lambda,p)P_{p-1,q} -a_{2}(\lambda,p,q-1)P_{p,q-2} \end{cases}
\eeqnas
which gives
\beqnas Z^{2}P_{p,q} &=& P_{p+2,q}
-\big(a_{1}(\lambda,p+1)+a_{1}(\lambda,p)\big)P_{p,q+1} \\&&
-\big(a_{2}(\lambda,p+1,q)+a_{2}(\lambda,p,q)\big) P_{p+1,q-1} \\
   && +a_{1}(\lambda,p)\big(a_{2}(\lambda,p-1,q+1)
    +a_{2}(\lambda,p,q)\big )  P_{p-1,q}\\&&
    +  a_{2}(\lambda,p,q)a_{2}(\lambda,p,q-1)P_{p,q-2}  \\
  &&+  a_{1}(\lambda,p)a_{1}(\lambda,p-1)P_{p-2,q+2}.\eeqnas
  On the other hand, from the induction hypothesis  we have
\beqnas && Z\Gamma(Z,P_{p,q})+ a_{1}(\lambda,p)\Gamma(Z,P_{p-1,q+1})+ a_{2}(\lambda,p,q)  \Gamma(Z,P_{p,q-1})=\\&&
\alpha_{0}(p,q)P_{p+2,q}+
\\&& \big[\alpha_{1}(p,q)-\alpha_{0}(p,q)a_{1}(\lambda,p+1) +\\&&
a_{1}(\lambda,p)\alpha_{0} (p-1,q+1)\big] P_{p,q+1}+\\
&&\big[\alpha_{2}(p,q)-\alpha_{0}(p,q)a_{2}(\lambda,p+1,q) +\\&&
a_{2}(\lambda,p,q)\alpha_{0} (p,q-1)\big] P_{p+1,q-1}+\\
&&\big[a_{1}(\lambda,p)\alpha_{2}(p-1,q+1)  +  \alpha_{1}(\lambda,p)  a_{2}(\lambda,p,q)- \\&&
 \alpha_{1}(p,q)  a_{2}(\lambda,p-1,q+1)- \alpha_{2}(p,q)  a_{1}(\lambda,p)\big]P_{p-1,q}+ \\
&&[a_{1}(\lambda,p)\alpha_{1}(p-1,q+1) - a_{1}(\lambda,p-1)\alpha_{1}(p,q)]P_{p-2,q+2}+\\
&&   \big[a_{2}(\lambda,p,q)\alpha_{2}(p,q-1) - a_{2}(\lambda,p,q-1)\alpha_{2}(p,q)\big] P_{p,q-2} , \eeqnas

Substituting everything in \eqref{recGamma}, we get
\beqnas \Gamma(Z,P_{p+1,q}) =&& (\alpha_{0}(p,q)-1) P_{p+2,q}+A_{1}(p,q) P_{p,q+1}+\\&&
A_{2}(p,q)P_{p+1,q-1}+A_{3}(p,q)P_{p-1,q}+\\&&
A_{4}(p,q)P_{p-2,q+2}+A_{5}(p,q)P_{p,q-2},\eeqnas
where
 \beqnas
 A_{1}(p,q)=&& 1+\alpha_{1}(p,q)-\alpha_{0}(p,q)a_{1}(\lambda,p+1) +\\&&a_{1}(\lambda,p)\alpha_{0} (p-1,q+1)\\&&+a_{1}(\lambda,p+1)+a_{1}(\lambda,p)\\
A_{2}(p,q) =&&-a_{1}(\lambda,q)+a_{2}(\lambda,p+1,q)+a_{2}(\lambda,p,q)+\alpha_{2}(p,q)\\&&-\alpha_{0}(p,q)a_{2}(\lambda,p+1,q) +a_{2}(\lambda,p,q)\alpha_{0} (p,q-1)\\
A_{3}(p,q) =&&-a_{2}(\lambda,q,p)-a_{1}(\lambda,p)a_{2}(\lambda,p-1,q+1) -a_{2}(\lambda,p,q)a_{1}(\lambda,p)  \\
&&+a_{1}(\lambda,p)\alpha_{2}(p-1,q+1)  +   \alpha_{1}(p,q-1)  a_{2}(\lambda,p,q)\\&& - \alpha_{1}(p,q)  a_{2}(\lambda,p-1,q+1)- \alpha_{2}(p,q)  a_{1}(\lambda,p)\\
A_{4}(p,q)=&& a_{1}(\lambda,p)\alpha_{1}(p-1,q+1) - a_{1}(\lambda,p-1)\alpha_{1}(p,q)\\&&
-  a_{1}(\lambda,p)a_{1}(\lambda,p-1)\\
A_{5}(p,q) =&& a_{2}(\lambda,p,q)\alpha_{2}(p,q-1) - a_{2}(\lambda,p,q-1)\alpha_{2}(p,q) \\
&&-a_{2}(\lambda,p,q)a_{2}(\lambda,p,q-1) .
 \eeqnas

A simple calculation shows that
\beqnas
&&1+\alpha_{1}(p,q)=  \alpha_{1}(p+1,q), \\
 &&A_{1}(p,q)=  \alpha_{1}(p+1,q), ~
 A_{2}(p,q)=  \alpha_{2}(p+1,q)\\
&& A_{3}(p,q)=  A_{4}(p,q)=  A_{5}(p,q)=0,\eeqnas
which concludes the induction formula for $\Gamma(Z, P_{p+1,q})$.
  The same method  leads to the formula for $\Gamma(\bZ,P_{p+1,q}) $, and exchanging $p$ and $q$ in the previous amounts to exchange $Z$ and $\bZ$.
\end{proof}
Now , we prove Theorem~\ref{defPpq}  using  Lemma~\ref{lemmaGammaPpq}.

\begin{proof}
Assume by induction  that $L^{(\alpha)}P_{p_1,q_1}=-\lambda_{p_1,q_1}P_{p_1,q_1} $   when  $p_1+q_1\leq p+q$, where $\lambda_{p,q}=(\lambda-1)(p+q)+p^2+q^2+pq$. As before,   we simply write the change of variable formula
\beqnas
\LL^{(\alpha)}P_{p+1,q}&=&\LL^{(\alpha)}(ZP_{p,q})+a_{1}(\lambda,p)\LL^{(\alpha)}P_{p-1,q+1}+a_2(\lambda,p,q)\LL^{(\alpha)}P_{p,q-1}\\
&=& Z\LL^{(\alpha)}P_{p,q}+P_{p,q}\LL^{(\alpha)}Z +2\Gamma(Z,P_{p,q})\\&&-a_{1}(\lambda,p)\lambda_{p-1,q+1}P_{p-1,q+1} -a_2(\lambda,p,q)\lambda_{p,q-1}P_{p,q-1}\\
&=& -(\lambda_{p,q}+\lambda)ZP_{p,q}+2\alpha_{0}(p,q)P_{p+1,q}+(2\alpha_{1}(p,q) \\&&-a_{1}(\lambda,p))P_{p-1,q+1} +(2\alpha_{2}(p,q)-a_{2}(\lambda,p,q)\lambda_{p,q-1})P_{p,q-1} .\eeqnas
But
$$Z P_{p,q}= P_{p+1, q} -a_{1}(\lambda,p)P_{p-1,q+1} -a_{2}(\lambda,p,q)P_{p,q-1},$$
so that
$$\LL^{(\alpha)}P_{p+1,q}=B_{1}(p,q)P_{p+1,q}+B_{2}(p,q)P_{p-1,q+1}+B_{3}(p,q)P_{p,q-1}, $$
where
\beqnas
B_{1}(p,q)&=& -(\lambda_{p,q}+\lambda)+2\alpha_{0}(p,q)\\
B_{2}(p,q)&=&2\alpha_{1}(p,q) +(\lambda_{p,q}-\lambda_{p-1,q+1}+\lambda)a_{1}(\lambda,p)\\
B_{3}(p,q)&=& 2\alpha_{2}(p,q)+(\lambda_{p,q}-\lambda_{p,q-1}+\lambda)a_{2}(\lambda,p,q)
\eeqnas Everything boils down to the following formulae, which are straightforward to check
$$ B_{1}(p,q)=-\lambda_{p+1,q},
 B_{2}(p,q)=B_{3}(p,q)=0.$$

The same proof applies for  $\LL^{(\alpha)}P_{p,q+1}=-\lambda_{p,q+1}P_{p,q+1}$. The conclusion follows.
\end{proof}
\brmq From the recurrence formula, it is easily checked that
\beqnas P_{p,q}=&& Z^{p}\bZ^{q}+A_{p,q}Z^{p+1}\bZ^{q-2}+B_{p,q}Z^{p-2}\bZ^{q+1}\\&&+C_{p,q}Z^{p-1}\bZ^{q-1}+D_{p,q}Z^{p-4}\bZ^{q+2}+F_{p,q}Z^{p+2}\bZ^{q-4} +R ,\eeqnas where ${\rm degree}(R)\leq p+q-3$. This general form may be easily induced from the form of the operator, and should produce a   six term recurrence formula. The fact  that the recurrence formula contains only 3 terms (as it is in dimension 1) is indeed quite mysterious.
\ermq

\section{Other representations of eigenpolynomials\label{sec.other} }
In this section, we come back to the two different representations for $\LL^{(-1/2)}$ and $\LL^{(1/2)}$ which provide new representations for the eigenvectors $P^{(\alpha)}_{p,q}$ in those specific cases.  This new representations will allow us to get in those cases linearization formulae for the product together with generating functions.

\subsection{Case $\alpha= -1/2$}
Although the $\alpha= -1/2$ case is quite easy,  since it comes from an Euclidean Laplace operator,  it gives rise to another family of recurrence formulae. On the other hand, the case $\alpha= 1/2$, which comes from the Casimir operator on $SU(3)$, leads to new representations of the eigenvectors $P^{(\alpha)}_{p,q}$ related to the irreducible representations of the symmetric group. In fine, comparing the two cases allows to  generalize the $SU(3)$ formulae to the general situation.

With the representation \eqref{matrix.complex2} of the operator $\LL^{(-1/2)}$, one may represent the function $Z$ as  a function $\R^2\mapsto \C$ as
$Z(x_1,x_2) = \frac{1}{3} (e_1+ e_j+ e_{\bar j})$, where, for $\beta= (\beta_1,\beta_2) \in \R^2$
$e_\beta= \exp\big(i(\beta_1x_1+\beta_2 x_2)\big)$, and $1, j, \bar j$ are the third root of unity, that is
$1= (1,0), ~j=(-1/2, \sqrt{3}/2), ~\bar j = (-1/2, -\sqrt{3}/2)$. Comparing with the description given in Section~\ref{sec.laplacien2d}, the change of normalization comes from the fact that we have divided $\LL$ by $4$ and replaced $Z$ by $Z/3$.

We have already observed that, for any triple $(z_1,z_2,z_3)$ of complex numbers satisfying $|z_i|=1$, $z_1z_2z_3=1$, setting for any $p\in \Z$,
$T_p= z_1^p+z_2^p+z_3^p$, one has  \eqref{recurrence.trace.su3}
\beq \label{recurrence.trace.su3} T_{p+2} -3ZT_{p+1} +3\bZ T_{p}-T_{p-1}= 0,\eeq
with $T_1= 3Z$ and $T_{-1}= 3\bZ$, $T_0=3$. One may observe first that this formula is unchanged if we replace $p$ by $-p$ and $Z$ by $\bZ$.
Setting $T_p= 3^{|p|} Q_p$, one gets
\beq\label{eqQp}Q_{p+1}= Z Q_p-\frac{1}{3} \bZ Q_{p-1} + \frac{1}{3^3} Q_{p-2}.\eeq
From this, it is clear that $Q_p$ is a polynomial with degree less than $p$ in the variables $(Z, \bZ)$, of the form $Q_p= Z^p+\hbox{\rm lower degree term}$.  Now, if we replace $z_1, z_2, z_3$ by $e_1, e_j, e_{\bar j}$ we see that $Q_p$ is an eigenvector for the Laplace operator $\Delta$ in $\R^2$, with eigenvalue $p^2$. Therefore,
\beq\label{eqQp} \forall p\geq 0, ~Q_p= P_{p,0}^{(-1/2)}, Q_{-p}=\overline{Q_p}=  P_{0,p}^{(-1/2)}.\eeq

Comparing~\eqref{eqQp} with the recurrence formulae for $P_{p,q}^{(\alpha)}$, the first line in formula~\eqref{defPpq2} gives in this case ($\lambda= 1$)
$$ P_{p+1,0}^{(-1/2)}= ZP_{p,0}^{(-1/2)}- \frac{1}{3} P^{(-1/2)}_{p-1,1},$$ which leads to
$$P_{p-1,1}^{(-1/2)}= \bZ P_{p-1,0}^{(-1/2)}-\frac{1}{3^2}P^{(-1/2)}_{p-2,0},$$ and the latter is nothing else than the second line in~\eqref{defPpq2}.

On the other hand, coming back to the representation $T_p= z_1^p+z_2^p+z_3^p$, one sees that, for any $(p,q)\in \Z^2$,
$$T_pT_q-T_{p+q}=\sum_{i\neq j} z_i^pz_j^q= \sum_{i\neq j} z_i^{p-q}z_j^{-q},$$ from which we get, for any $(p,q)\in \Z^2$
\beq\label{eqTpq}T_pT_q-T_{p+q}=  T_{p-q}T_{-q}-T_{p-2q}= T_{q-p}T_{-p}-T_{q-2p}.\eeq
 When $(z_1,z_2,z_3)= ( e_1, e_j,e_{\bar j})$, this writes as sums of terms of the form $e_{\beta(p,q)}$, where $|\beta_{p,q}|^2= p^2+q^2-pq$. Therefore, for the $\LL^{(-1/2)}$ operator, writing $T_pT_{-q}-T_{p-q}$ as a polynomial in $(Z, \bZ)$, we see that  this is an eigenvector associated with the eigenvalue $p^2+q^2+pq$. Looking at the highest degree term, and translating this in terms of the polynomials $Q_p= 3^{-|p|} T_p$, we obtain

\beq\label{eqvpplat}\forall p, q\geq 1, ~P^{(-1/2)}_{p,q}=Q_pQ_{-q}- 3^{-2\min(p,q)}Q_{p-q},\eeq
which gives a representation  of $P_{(p,q)}^{(-1/2)}$ in terms of the polynomials $P_{p,0}^{(-1/2)}$ and $P_{0,q}^{(-1/2)}$ which is not easy to obtain directly from the recurrence formula~\eqref{defPpq2}.

When $p,q\geq 0 ,Q_pQ_q-Q_{p+q}$  is  also an eigenvector for the Laplace operator associated with the eigenvalue $p^2+q^2-pq $. Indeed, using~\eqref{eqTpq}
which is valid for any $(p,q)\in \Z^2$, and  comparing with~\eqref{defPpq2},
we end up, for $p\ge q\ge 0$,  with
$$ P^{(-1/2)}_{p-q,q} =Q_pQ_q-Q_{p+q},$$

\bprop[Linearization formula]

\beqnas
P^{(-\frac{1}{2})}_{p,q}P^{(-\frac{1}{2})}_{p',q'}=&&P^{(-\frac{1}{2})}_{p+p',q+q'}+\\
&&3^{-\min(q,q')}P^{(-\frac{1}{2})}_{p+p'+\min(q,q'),\max(q,q')-\min(q,q')}\\&&+3^{-\min(p,p')}P^{(-\frac{1}{2})}_{\max(p,p')-\min(p,p'),q+q'+\min(p,p')}\\&&
+b_{1}(p,q,p',q')P^{(-\frac{1}{2}}_{|\gamma|-\max(0,\delta),|\delta|+\min(0,\gamma)}\\&&+b_{2}(p,q,p',q')P^{(-\frac{1}{2})}_{|\beta|+min(0,\alpha),|\alpha|+min(0,\beta)}\\&&
+b_{3}(p,q,p',q')P^{(-\frac{1}{2})}_{|\beta'|+min(0,\alpha'),|\alpha'|+min(0,\beta')}
\eeqnas
Where
\beqnas
\gamma&=&\max(p'-q,p-q')\\
\delta&=&\min(p'-q,p-q')\\
\alpha&=&p+q-q'\\
\beta&=&p'+q'-p\\
\alpha'&=&p'+q'-q\\
\beta'&=&p+q-p'\\
b_{1}(p,q,p',q')&=&3^{|\gamma|+|\delta|+\min(0,\gamma)-\max(0,\delta)-(p+q+p'+q')}\\
b_{2}(p,q,p',q')&=&3^{|\beta|+|\alpha|+min(0,\alpha)+min(0,\beta)-(p+p'+q+q')}\\
b_{3}(p,q,p',q')&=&3^{|\beta'|+|\alpha'|+min(0,\alpha')+min(0,\beta')-(p+p'+q+q')}\\
\eeqnas
\eprop

\subsection{Case $\alpha= 1/2$.}
We now turn to the inspection of the family $P_{p,q}^{(1/2)}$. We know that $\frac{4}{3}\LL^{(1/2)}$  may be represented as the action of the Casimir operator on $SU(3)$ acting on spectral functions.
Comparing with formulae~\eqref{eqLTpSU3} and ~\eqref{eqGammaTpSU3}, one sees that, for $p\geq 1$,
\beq\label{eqLTpL12}\LL^{(1/2)} (T_p)= -\frac{p}{4}\Big(2(9-p)T_p+ 3\sum_{i=1}^{p-1} T_iT_{p-i}\Big),\eeq
while, for $p,q\geq 0$
\beq\label{eqGammaTp} \Gamma( T_p,T_q)= \frac{pq}{2}\big(T_pT_q-3T_{p+q}\big), \eeq
with similar formulae for $p\leq 0, q\leq 0$.

If we remember that
$\LL^{(1/2)}= \LL^{(-1/2)}-3(Z\partial_Z+ \bZ\partial_{\bZ})$, we end up with the formula for $p\geq 0$

\beq\label{eqZTp}\big(Z\partial_Z+ \bZ\partial_{\bZ}\big)T_p= -\frac{p}{2}(p-3)T_p+ \frac{p}{4}\sum_{i=1}^{p-1} T_iT_{p-i}.\eeq

With the help of formula \eqref{comp.Lalpha},  for $p\geq 1$ we end up with
\beq\label{eqLalphaTp} \LL^{(\alpha)}T_p = -\frac{p}{4}\big(p(1-6\alpha)+ 9(2\alpha+1)\big)T_p-\frac{3p}{8}(2\alpha +1)\sum_{i=1}^{p-1} T_iT_{p-i}.\eeq

Fix now some integer $n$ and denote by $\Pi_n$ the set of sequences   $\pi= (p_1, \cdots, p_k)$ of integers $ p_1\geq p_k\geq 1$  such that $p_1+ \cdots + p_k= n$. For $\pi\in \Pi_n$,  denote $T_\pi= \prod_{j=1}^k T_j$.
 Comparing with formula~\eqref{eqGammaTp}, and the general formula~\eqref{formule.produit},  we see that the vector space generated by $T_\pi, \pi\in \Pi_n$, is stable under $\LL^{(\alpha)}$. We therefore will be able to diagonalize $\LL^{(\alpha)}$ in this vector space.

 We first perform a slight change in the normalization of the variables $T_p$, setting $T_p= c S_p$, in order to reduce formulae ~\eqref{eqLalphaTp} and~\eqref{eqGammaTp} to
 \beq \LL^{(\alpha)}S_p = -\mu_{p, \alpha} S_p -\frac{1}{c} \frac{3p}{4} \sum_{i=1}^{p-1} S_i S_{p-i},\eeq
   and
 \beq \Gamma(S_p,S_q)= \frac{pq}{2} (S_pS_q-\frac{3}{c} S_{p+q}).\eeq

 With $\mu_{p, \alpha} = \frac{p}{4}\big(p(1-6\alpha)+ 9(2\alpha+1)\big)$ and   $c=\sqrt{\frac{2}{2\alpha+1}}$,

Following \cite{T.Levy},  it is easier to introduce the group $\cS_n$ of order $n$ permutations. For any $\sigma\in \cS_n$, one consider it's cycle decomposition $\sigma= \sigma_1\cdots \sigma_k$ (ordered in increasing lengths) and denote by $\pi= \pi(\sigma)= (p_1, \cdots, p_k)\in \Pi_n$ the sequence of the lengths of $\sigma_j$. We then denote  $S_\sigma= S_{p_1}\cdots S_{p_k}$. It is worth to observe that if $\tau= (ij)$ is a transposition, the cycle decomposition of $\sigma\tau$ splits one cycle in two subcycles when $i$ and $j$ belong to the same cycle and glues together two cycles when $i$ and $j$ belong to two different cycles.

Therefore, if $\cT_n$ denotes the set of all transpositions, for a permutation $\sigma$ with  $\pi(\sigma)=(p_1, \cdots, p_k)$,  through an easy combinatorial argument, one gets
$$\sum_{\tau\in \cT_n} S_{\tau\sigma} = \frac{1}{2}\sum_{i=1}^k p_i\frac{S_\sigma}{S_{p_i}}\sum_{j=1}^{p_i-1} S_jS_{p_i-j}+ \sum_{i,j=1}^kp_ip_j\frac{S_\sigma}{S_{p_i}S_{p_j}} S_{p_i+p_j}.$$

Comparing this and the formula~\eqref{formule.produit}, we get
\beq \LL^{(\alpha)}(S_\sigma)= -\mu_{\sigma, \alpha} S_\sigma -\frac{3}{2c} \sum_{\tau\in \cT_n} S_{\sigma\tau},
\eeq

More precisely
\beqnas
\LL^{(\alpha)}S_{\sigma}&=&
S_{\sigma}\Big(-\sum_{i=1}^{k}\mu_{p_i, \alpha}+\sum_{i \neq j}\frac{p_ip_j}{2}\\&&-\frac{3}{2c}
\big( \sum_{i=1}^{k} \frac{p_i}{2S_{p_i}}\sum_{j=1}^{p_i-1}S_jS_{p_i-j} +\sum_{i\neq j} p_ip_j \frac{S_{p_i+p_j}}{S_{pi}S_{p_j}}\big)   \Big)\\
&=&
 S_{\sigma}\Big(-\mu_{\sigma, \alpha}-\frac{3}{2c}
\big( \sum_{i=1}^{k} \frac{p_i}{2S_{p_i}}\sum_{j=1}^{p_i-1}S_jS_{p_i-j} +\sum_{i\neq j} p_ip_j \frac{S_{p_i+p_j}}{S_{pi}S_{p_j}}\big)   \Big)
\eeqnas

 where, for  $\pi(\sigma)= (p_1, \cdots, p_k)$,
\beqnas
\mu_{\sigma,\alpha}&=&\sum_{i=1}^{k}\mu_{p_i, \alpha}-\sum_{i \neq j}\frac{p_ip_j}{2}\\
&=&\frac{3}{4}(1-2\alpha)\sum_{i=1}^{k}p^{2}_i+\frac{9}{4}(2\alpha +1)n-\frac{n^2}{2}
\eeqnas
Finally
\beqnas
\LL^{(\alpha)}S_{\sigma}&=&-\mu_{\alpha,\sigma}S_{\sigma}-\frac{3}{2c}\sum_{\tau\in \cT_n} S_{\tau\sigma}
\eeqnas

where $\cT_n$ is the set of transpositions in $\cS_n$.

It is worth to observe that for $\alpha= 1/2$ (and only in this case),  $\mu_{\sigma, \alpha}$ depends on $n$ only, and therefore finding eigenvectors for $\LL^{(1/2)}$ amounts to find eigenvectors for the linear operator
$S_{\sigma}\mapsto \sum_{\tau\in \cT_n} S_{\sigma\tau}$. But the latter corresponds to the  operator $\sum_{\tau\in \cT_n}\tau$ in the group algebra of the group $\cS_n$, which commutes to every group element. It is therefore  diagonal on any irreducible representation. Turning back to our setting, we conclude that for any character $\chi$ of the group $\cS_n$,
$\sum_{\sigma\in \cS_n} \chi(\sigma) S_\sigma$ is an eigenvector for $\LL^{(1/2)}$. The characters of the group $\cS_n$ are well known and correspond to Young diagrams, which are indeed elements of the set $\Pi_n$ described above.  Unfortunately, as we shall see in the following examples, this representation is far from being one to one, and many eigenvectors coming from the $\cS_n$ representation has degree less than $n$ in $(Z, \bar Z)$. The correspondence between the degree and the shape of the Young diagram remains quite mysterious.

 The paper \cite{Caracter} describes an elegant method  which provides  a simple combinatorial way for  computing the character table in  any symmetric group $S_n$. Since for any character $\chi$, the value of $\chi(\sigma)$ only depends on the conjugacy class of $\sigma$,  that is on the Young diagram it belongs to, one has to compute $\chi(\xi)$ for any pair $(\chi, \xi)$ of Young diagrams. This may be achieved through the analysis of the so-called border strips.  In what follows, we then give some examples of eigenvectors for $\LL^{(1/2)}$ provided by this description, that is $\sum_{\pi\in \Pi_n} |\pi|\xi(\pi)  T_\pi$, where $|\pi|$ denotes the size of the conjugacy class $\pi \in \cS_n$, that is the number of those elements $\sigma\in \cS_n$  such that $\pi(\sigma)= \pi$. For $\pi= (p_1, \cdots, p_k)$, one has
 $|\pi|= \dis \frac{n!}{\prod_{j=1}^{n}k_j!}\prod_1^k\frac{1}{p_j}$,  where $k_j$ is the number of cycles with length $j$ in $\pi$.

 As an example, we show below the eigenvectors given by this construction for $n=2,3,4$

The group $\cS_2$ has two conjugacy classes $\chi_1, \chi_2$, corresponding to the partitions ($2,0)$ and $(1,1)$,
$$
\chi_1=\begin{Young}
& \cr
\end{Young}\,,~ \chi_2=
\begin{Young}
\cr
\cr
\end{Young}\,.
$$
The character table is then
$$
\begin{tabular}{|c|cc|}
    \hline
     & $\chi_1$ & $\chi_2$   \\
    \hline
    $\chi_1$ & 1 & 1  \\
    $\chi_2$ & 1 & -1  \\
    \hline
\end{tabular}
$$
\\
From his we get the following eigenvectors for  $\LL^{(1/2)}$:
\beqnas
Q_{1}(Z,\bZ)&=&\chi_{1}(\chi_1)T_{2}+\chi_{1}(\chi_2)T_1^{2}\\
&=&T_{2}+T^{2}_1\\
&=&18Z^2-6\bZ
\eeqnas

\beqnas
Q_{2}(Z,\bZ)&=&\chi_{2}(\chi_1)T_{2}+\chi_{2}(\chi_2)T_1^{2}\\
&=&T^{2}_1-T_{2}\\
&=&6\bZ
\eeqnas
For $\cS_3$, we have   three conjugacy classes   $\chi_1, \chi_2, \chi_3$ corresponding to the partitions  $(3,0,0), (2,1,0), (1,1,1)$

$$
\chi_1=\begin{Young}
 &&\cr
\end{Young}\,, ~
\chi_2=\begin{Young}
 \cr \\
 &\cr
\end{Young}\,, ~
\chi_3=\begin{Young}
 \cr\\
 \cr\\
 \cr
\end{Young}\,,
$$
with character table
$$
\begin{tabular}{|c| ccc| }
    \hline
     & $\chi_1$ & $\chi_2$&$\chi_3$  \\
    \hline
    $\chi_1$ & 1 & 1 & 1  \\
    $\chi_2$ & 2 & 0 & -1  \\
    $\chi_3$ & 1 & -1& 1\\
    \hline
\end{tabular}
$$
and corresponding eigenvectors
\beqnas
Q_1(Z,\bZ)&=& \chi_{1}(\chi_1)T_1^{3}+\chi_{1}(\chi_2)T_2T_1+\chi_{1}(\chi_3)T_{3}\\
&=&T_1^{3}+T_2T_1+T_3\\
&=& 3(27Z^3-15Z\bZ+1)
\eeqnas

\beqnas
Q_2(Z,\bZ)&=& \chi_{2}(\chi_1)T_3+\chi_{2}(\chi_2)T_1T_2+\chi_{2}(\chi_3)T^2_{1}\\
&=&2T_3-2T_1^{3}\\
&=& 2(27Z\bZ-3)
\eeqnas
\beqnas
Q_3(Z,\bZ)&=& \chi_{3}(\chi_1)T_3+\chi_{3}(\chi_2)T_1T_2+\chi_{3}(\chi_3)T^2_{1}\\
&=&T_1^{3}-3T_2T_1+2T_3\\
&=& 6
\eeqnas
For  $\cS_4$, we  have five conjugacy classes,   with Young diagrams $\chi_1, \chi_2, \chi_3, \chi_4, \chi_5$,  corresponding to $(4,0,0,0), (2,1,1,0), (2,2,0,0), (3,1,0,0), (1,1,1,1)$.
 $$
\chi_1=\begin{Young}
 &&&\cr
\end{Young}\,,~
\chi_2=\begin{Young}
\cr\\
\cr\\
&\cr
\end{Young}\,,~
\chi_3=\begin{Young}
 &\cr \\
 &\cr
\end{Young}\,, ~
\chi_4=\begin{Young}
 \cr \\
 &&\cr
\end{Young}\,, ~
\chi_5=\begin{Young}
 \cr \\
 \cr\\
 \cr\\
 \cr
\end{Young}\,,
$$
and character table
$$
\begin{tabular}{|c| ccccc| }
    \hline
     & $\chi_1$ & $\chi_2$&$\chi_3$ &$\chi_4$&$\chi_5$ \\
    \hline
    $\chi_1$ & 1 & 1 & 1 & 1& 1 \\
    $\chi_2$ & 1 & -1 & 1 & 1 & -1 \\
    $\chi_3$ & 2 & 0 & 2 & -1 & 0 \\
   $\chi_4$ & 3 & 1 & -1 & 0 &-1\\
   $\chi_5$ & 3 & -1 & -1 & 0 & 1\\
    \hline
\end{tabular}
$$
The  corresponding eigenvectors are then:
\beqnas
Q_1(Z,\bZ)&=& \chi_1(\chi_1)T^{4}_1+ 6 \chi_1(\chi_2)T^2_1T_2+3 \chi_1(\chi_3)T^2_2+ 8\chi_1(\chi_4)T_1T_3\\
 &&+6\chi_1(\chi_5)T_4\\
&=& T^{4}_1+ 6 T^2_1T_2+ 3 T^2_2+ 8T_1T_3+ 6T_4\\
&=&72(27Z^4-27\bZ Z^2+3\bZ^2+2Z)
\eeqnas

\beqnas
Q_2(Z,\bZ)&=& \chi_2(\chi_1)T^{4}_1+6 \chi_2(\chi_2)T^2_1T_2+ 3 \chi_2(\chi_3)T^2_2+ 8\chi_2(\chi_4)T_1T_3\\
 &&+6\chi_2(\chi_5)T_4\\
&=& T^{4}_1- 6 T^2_1T_2+3 T^2_2+ 8T_1T_3- 6T_4\\
&=&0
\eeqnas
\beqnas
Q_3(Z,\bZ)&=& \chi_3(\chi_1)T^{4}_1+ 6 \chi_3(\chi_2)T^2_1T_2+ 3 \chi_3(\chi_3)T^2_2+8\chi_3(\chi_4)T_1T_3\\
&& +6\chi_3(\chi_5)T_4\\
&=& 2T^{4}_1+2.3 T^2_2-8T_1T_3\\
&=&72(3\bZ^2-Z)
\eeqnas
\beqnas
Q_4 (Z,\bZ)&=& \chi_4(\chi_1)T^{4}_1+ 6 \chi_4(\chi_2)T^2_1T_2+3 \chi_4(\chi_3)T^2_2+ 8\chi_4(\chi_4)T_1T_3\\
&&+ 6\chi^4(\chi_5)T_4\\
&=& 3T^{4}_1+6 T^2_1T_2- 3 T^2_2- 6T_4\\
&=&72(9\bZ Z^2-3\bZ^2-Z)
\eeqnas
\beqnas
Q_5(Z,\bZ) &=& \chi_5(\chi_1)T^{4}_1+6 \chi_5(\chi_2)T^2_1T_2+ 3 \chi_5(\chi_3)T^2_2+{8}\chi_5(\chi_4)T_1T_3\\
 &&+6\chi_5(\chi_5)T_4\\
&=& 3T^{4}_1- 6 T^2_1T_2- 3 T^2_2+ 6T_4\\
&=&72 Z
\eeqnas

\section{Generating functions }

In this section, we provide first a partial generating function in the general case for the family $P^{(\alpha)}_{0,n}$ or equivalently for $P^{(\alpha)}_{n,0}$ , which leads to some simple representation of the polynomials $P_{n,m}^{(\alpha)}$ as linear combinations of $P^{(\alpha)} _{p,0}P^{(\alpha)}_{0,q}$. then we turn to the two geometric cases in which a complete generating function may be provided.

\subsection{Partial generating functions in the general case\label{subsec.generating.gal}}
In this section, we propose an alternative representation of the eigenvectors in the general case, together with a partial generating function.   We start with the following

\bprop\label{prop.LPalpha}

Let
\beq
P(X)= 1-3\bZ X+3Z X^2-X^3 ,
\eeq

Then, still with $\lambda=\frac{1}{2}(6\alpha+5)$,
\beqnas
\LL^{({\alpha})}(P(X))&=&-\lambda XP'+\frac{\lambda}{2}X^2P''
  ,\\
  \Gamma(P(X),P(Y))&=&\frac{XY}{2}\big(P'(X)P'(Y)+3\frac{P'(X)P(Y)-P(X)P'(Y)}{X-Y}\big),\\
  \hbox{ from which}\\
 \Gamma(P(X),P(X))&=& \frac{X^2}{2}(3PP''-2P'^2),\\
\eeqnas

Also, with $\bar P(Y)= 1-3Z Y+3 \bZ Y^2-Y^3=-Y^3P(1/Y)$
$$\Gamma(P(X), \bar P(Y))= \frac{XY}{2(XY-1)}\Big(3XP'\bar P+3Y\bar P'P-9P\bar P-(XY-1)P'\bar P'\Big).$$
\eprop
\bpf
The proof boils down to a simple verification,  using  the linearity of $\LL$ and the bilinearity of $\Gamma $.
The formula for  $\Gamma(P(X),P(X))$ may be obtained directly from $\Gamma(P(X),P(X))=\lim_{Y \to X}  \Gamma(P(X),P(Y)).$
\epf
 \bprop\label{prop.fct.gen2}

 Let $Q=P^{\beta}$, $\beta=\frac{1-\lambda}{3}= -\frac{1+2\alpha }{2}$.  Then, for $\alpha\neq -1/2$, one has

 \beq
\label{eqLaQ} \LL^{(\alpha)}(Q(X))=-\lambda XQ'(X)-X^2Q''(X),\eeq

\beq \label{eqGammaaQ} \Gamma(Q(X), Q(Y))=\frac{XY}{2}\big(Q'(X)Q'(Y)
+3\beta\frac{Q'(X)Q(Y)-Q(X)Q'(Y)}{X-Y}\big).
  \eeq
   \beq\label{eqGammaaQb}
 \Gamma(Q(X), \bar Q(Y)) =3\beta^{2} \frac{XY}{2(XY-1)}Q(X)\bar Q(Y)\Big(XS+Y\bar{S}- 3\Big) -\frac{XY}{2}Q'(X)\bar{Q}'(Y),
    \eeq\\
    where $S=\frac{P'}{P}(X),~~ \bar{S}=\frac{\bar P'}{\bar P}(Y)$\\
 \Big( For $\alpha= -1/2$, one should  replace $ Q= P^\beta$ by $Q=\log(P)$\Big).
 \eprop

 \bpf
 Let us look first at $\LL^{(\alpha)}(Q)$.
 With the formulae \eqref{change.var}, we have
 \beqnas\label{derive de Q}
 \LL^{(\alpha)}(Q(X))&=&\beta P(X)^{\beta-1}  \LL^{(\alpha)}(P(X))+\beta(\beta-1)P(X)^{\beta-2} \Gamma(P(X),P(X)),
 \eeqnas
 and from Proposition~\ref{prop.LPalpha}, we have
 \beqnas
 \LL^{(\alpha)}(Q(X))&=& -\lambda X Q'(X) +X^2 [(\frac{\lambda}{2}+\frac{3}{2}(\beta-1))\beta P''(X)P(X)^{\beta-1}\\
 &&-\beta (\beta-1)P'^{2}(X)P^{\beta-2}(X)].
 \eeqnas
For  the particular value of $\beta=\frac{1-\lambda}{3}$, $\frac{\lambda}{2}+\frac{3}{2}(\beta-1)=-1$, and  we get  the announced result.

 Turning now to formula~\eqref{eqGammaaQ}, we write
 \beqnas
 \Gamma(Q(X),Q(Y))&=& \beta^2P^{\beta-1}(X)P^{\beta-1}(Y)\Gamma(P(X),P(Y))\\
 &=&\frac{XY}{2}\Big(   \beta^{2}P'(X)P^{\beta-1}(X)P'(Y)P^{\beta-1}(Y)\\&&
  +\frac{3\beta}{X-Y}\big(\beta P^{\beta}(Y)P'(X)P^{\beta-1}(X)\\
 && -\beta P^{\beta}(X)P'(Y)P^{\beta-1}(Y)\big)\Big)
 \eeqnas
 It remains to write
  $ Q'=\beta P'P^{\beta-1}$ to obtain  formula~\eqref{eqGammaaQ}.

  Formula~\eqref{eqGammaaQb}   is obtained in the same way.

 \epf

\bcor  \label{cor.gen.gal} $Q(X)$ is a generating function for the family $P_{0,n}^{(\alpha)}$. More precisely, still with $\beta= -(1+2\alpha)/2$,
$$(1-3\bar Z X+ 3Z X^2-X^3)^\beta= \sum_n c_nP_{0,n}^{(\alpha)} X^n,$$ where
$$c_n={\displaystyle (-3)^{n}\frac{\beta(\beta-1)....(\beta-n+1 )}{n!}}$$\\
\big(The same remarks as in Proposition~\ref{prop.fct.gen2} applies for $\alpha= -1/2$\big).
\ecor

\bpf
If we write the asymptotic expansion of $Q(X)$ around $X=0$ (which is licit since both $Z$ and $\bar Z$ are bounded), and writing
$Q(X)= \sum_n  A_n(Z, \bar Z) X^n$, equation~\eqref{eqLaQ} gives
$$\LL ^{(\alpha)}A_{n}=-(\lambda n+2n(n-1))A_{n}.$$

As a consequence,  then $A_n(Z,\bZ)$ are  eigenvectors of $\LL^{(\alpha)}$. But a simple computation shows that $A_n(Z, \bar Z)$ is a polynomial in $(Z, \bar Z)$ with   highest degree term $c_n \bar Z^n$. Therefore, $A_n= c_n P_{0,n}^{(\alpha)}$.

\epf
For example:
\beqnas
A_{0}(Z,\bZ)&=&1\\
A_{1}(Z,\bZ)&=&-3 \beta \bZ\\
A_{2}(Z,\bZ)&=&3\beta(3(\beta-1)\bZ^2+2Z)
\eeqnas

If we set  $\bar Q= \bar P^\beta$, then $Q(X)\bar Q(Y)= \sum_{p,q} c_pc_q P_{0,p}^{(\alpha)} P_{q,0}^{(\alpha)} X^pY^q$, where the highest degree term in the polynomial $P_{0,p}^{(\alpha)}P_{q,0}^{(\alpha)}$ is $Z^q\bar Z^p$. Considering similarly a differential equation satisfied by $Q(X)\bar Q(Y)$ will then provide useful informations about the eigenvectors $P^{(\alpha)}_{p,q}$. \\

Obtaining  two variables generating functions for the general family $P_{m,n}^{(\alpha)}$ is not easy.
To simplify the calculations, it is simpler to  introduce the operator
$\hat L= \LL^{(\alpha)}+L_{0}$, where
$$ L_{0}=\lambda(X\partial_{X}+Y\partial_{Y})+X^2\partial_{X}^2+Y^2\partial_{Y}^2+XY\partial_{X,Y}^2,$$ and associated  $\Gamma_0$ operator
 $$\Gamma_{0}(f,g)= X^2\partial_{X}f \partial_{X}g +Y^2\partial_{Y}f \partial_{Y}g + \frac{1}{2}XY(\partial_{X}f\partial_{Y}g+\partial_{X}g\partial_{Y}f)
$$
and $\hat \Gamma= \Gamma+ \Gamma_0$.

Then, from Proposition~\ref{prop.fct.gen2}, we have $\hat \LL(Q(X))= \hat \LL(\bar Q(Y))=0$ and therefore
\beqnas
\hat L (Q(X)\bar Q(Y))&= &2\hat \Gamma(Q(X),\bar Q(Y))\\
&=& 2\big(\Gamma(Q(X),\bar Q(Y))+\frac{XY}{2} Q'(X) \bar Q'(Y)\big),
\eeqnas

form which we get
\bprop\label{form.hLQbQ}
$$
\hat L (Q(X)\bar Q(Y))=\frac{3\beta^{2} XY}{XY-1}Q(X)\bar Q(Y)\Big( XS+Y\bar S  -3 \Big),
$$
where as before $S=\frac{P'}{P}(X)$, $\bar S=\frac{\bar P'}{\bar P}(Y).$
\eprop

This proposition leads us to a new representation of the polynomials $P^{(\alpha)}_{n,m}$.
Writing for simplicity $Q(X)=\sum_{n\ge0}A_nX^n$ and  $\overline{Q}(Y)=\sum_{m\ge0}B_mY^m$ then we deduce that :
\beqnas
\LL^{(\alpha)}(A_nB_m-A_{n-1}B_{m-1})&=-&\lambda_{n,m}A_{n}B_{m}-\delta_{n,m}A_{n-1}B_{m-1}
\eeqnas
where $$\lambda_{n,m}=(\lambda-1)(n+m)+n^2+m^2+nm.$$

	and
	$$ \delta_{n,m}=(1-\lambda)(\lambda+n+m-3)-\lambda_{n-1,m-1}.$$

Finally,  from an easy induction, one sees that  $L^{(\alpha)}(A_nB_m)$ is a linear combination of  $A_{n-p}B_{m-p}$, $0\leq p \leq \min(n,m)$, from which one deduces that $P_{m,n}^{(\alpha)}$ may be written as
$$P_{m,n}^{(\alpha)}= \sum_{p=0}^{\min(m,n)} d_{m,n,p, \alpha} P_{m-p,0}^{(\alpha)}P_{0,n-p}^{(\alpha)}$$
very similar to the representation given for the case $\alpha= -1/2$ in equation~\eqref{eqvpplat}. Unfortunately, the explicit expression for the constants $d_{m,n, p, \alpha}$ does seem to have any simple form.

 From what preceedes,  we see that  the family $Q(X)\bar Q(Y)$  is not a generating function for  the family
 $P_{m,n}^{(\alpha)}(Z,\bZ)$, but we may expect that it is the case for    some expression of the form $F(XY)Q(X)\bar Q(Y)$ for some real valued function $F$. Therefore, we may look at  an equation of the form
 $\hat{L}(F(XY)Q\bar{Q})=0$. We may then use the following remark
 \bprop
\beqnas
\hat {\Gamma}(XY, Q(X)\bar Q(Y))=\frac{3}{2}\beta XY Q\bar Q(XS+Y\bar S)\\
\eeqnas
where $S=\frac{P'}{P}(X)$, $\bar S=\frac{\bar{P}'}{\bar{P}}(Y)$
\eprop
\bpf
\beqnas
\hat{\Gamma}(XY, Q(X)\bar Q(Y))&=& \Gamma_{0}(XY,Q\bar Q)\\
&=&X\big(Q\Gamma_{0}(Y,\bar Q)+\bar Q\Gamma_{0}(Y,Q)\big)+\\&&Y\big(  Q\Gamma_{0}(X,\bar Q)+\bar Q\Gamma_{0}(X,Q)\big)\\
&=&\frac{3}{2}XY\big(YQ(X)\bar Q'(Y) +X\bar Q(Y)Q'(X)  \big)\\
&=&\frac{3}{2}\beta XY Q(X)\bar Q(Y) \Big(XS+Y\bar S\Big)
\eeqnas

\epf
Then we have,
\beqnas
\hat L(F(XY)Q\bar Q)&=&F(XY)\hat L(Q\bar Q)+Q\bar Q \hat L(F(XY))+2\hat{\Gamma}(F(XY), Q\bar Q)\\
&=& F(XY)\hat L(Q\bar Q)+Q\bar Q\big( F'(XY)\hat{L}(XY)+\\&&F''(XY)\hat{\Gamma}(XY,XY)\big)
 +2F'(XY)\hat{\Gamma}(XY,Q\bar Q)
\eeqnas
A simple computation  shows that $\hat{L}(XY)=(2\lambda+1)XY$, $\hat{\Gamma}(XY,XY)=3X^2Y^2$.
Finally we have
\beqnas
\hat L(F(XY)Q\bar Q)&=&2\hat{\Gamma}(F(XY),Q\bar{Q})\big(F'(XY)+\frac{\beta}{XY-1}F(XY)\big)\\&&
+ Q\bar{Q}XY\Big(  (2\lambda+1)F'(XY)+3XYF''(XY)-\\
&&\frac{9\beta^2}{XY-1}F(XY)\Big)
\eeqnas
In order to get  $\hat L(F(XY)Q\bar Q)=0$, two  are led to solve the two linear differential equations on $F$
\beqnas
   \begin{cases}
   F'(u)+\frac{\beta}{u-1}F(u)= 0\\
   \\
   3uF''(u)+(2\lambda+1)F'(u)-\frac{9\beta^2}{u-1}F(u)=0\\
   \end{cases}
   \eeqnas
and it is easy to see that they are compatible only when
$\lambda^2-5\lambda+4=0$ that is  $\lambda \in \{1,4\}$, which leads to the study of the two geometric cases.
 \subsection{Case ${\alpha=-1/2}$}

 In this case, the previous approach just provides $\log(P(X))-\log (\bar P(Y))$ as a bivariate  generating function, which is too degenerate to give useful information on the general polynomials $P_{n,m}^{(-1/2)}$. But the explicit representation of the eigenvectors in this case pearly provides an efficient generating function already obtained in {\cite{dunn}}.

\bprop
A generating function in the case $\alpha= -1/2$  for the family $S^{(-\frac{1}{2})}_{p,q}= T_pT_{-q}-T_{p-q}$,  ~~  $p,q \geq 0$ is defined by :
\beqnas
G(X,Y)&=&\big(3-X\frac{\bar P'}{\bar P}(X)\big) \big( 3-Y\frac{ P'}{ P}(Y)\big)\\
&&+\frac{1}{1-XY}\Big(  X\frac{\bar P'}{\bar P}(X)+ Y\frac{P'}{ P}(Y)-3 \Big)
 \eeqnas
  \eprop
 \bpf
  It is quite easy to deduce the generating function  for the family $S^{(-\frac{1}{2})}_{p,q}$. Indeed, since $S^{(-\frac{1}{2})}_{p,q}= T_pT_{-q}-T_{p-q}$,  ~~  $p,q \geq 0$  we get
\beqnas G(X,Y)&:=&\sum_{p\geq 0, q\geq 0} X^pY^q S_{p,q}\\&=& (\sum_{p\geq 0} T_pX^p)(\sum_{q\geq 0} T_{-q}Y^q)- \sum_{p\geq 0,p\geq 0}T_{p-q} X^{p}Y^{q},\eeqnas the series being convergent as soon as $|X|<1$ and $|Y|<1$. Using the representation $T_p = z_1^p+z_2^p+z_3^p$, this sums as
$$\Big(\sum_{i=1}^3(1-Xz_i)^{-1}\Big)\Big(\sum_{j=1}^3(1-Y\overline{z_i})^{-1})\Big)-\sum_{i=1}^3(1-Xz_i)^{-1}(1-Y\overline{z_i})^{-1},$$
But
$$\sum_{i=1}^3(1-Xz_i)^{-1}= \frac{1}{X} \frac{P'}{P}(\frac{1}{X}).$$
On the other hand
$$ (1-Xz_i)^{-1}(1-Y\overline{z_i})^{-1}= \frac{1}{1-XY}(\frac{1}{1-Xz_i}+ \frac{1}{1-Y\bar z_i}-1),$$
so that
$$\sum_{i=1}^3(1-Xz_i)^{-1}(1-Y\overline{z_i})^{-1}= \frac{1}{1-XY}\Big( \frac{1}{X} \frac{P'}{P}(\frac{1}{X})+  \frac{1}{Y}\frac{\bar P'}{\bar P}(\frac{1}{Y}) -3\Big).$$
This simplifies using $X^3P(1/X)= -\bar P(X)$, such that
$$ \frac{1}{X} \frac{P'}{P}(\frac{1}{X})= 3-X\frac{\bar P'}{\bar P} (X),$$

 Then we have the result. Using the notations of the previous subsection~\ref{subsec.generating.gal}, we may also check directly that
 $(\LL^{(-\frac{1}{2})}+L_{0})G(X,Y)=0$ in this case.

 \epf

 \subsection  {Case $ \alpha=1/2$}

Finally, the formulae providing in subsection~\ref{subsec.generating.gal} leads us directly to a generating function in the $\alpha= -1/2$ case. This generating function has been proposed in \cite {dunn} with however a completely different approach, based on the representations of the $SU(3)$ group.
 \bprop
A generating function for the family $P^{(\frac{1}{2})}_{m,n}$ is given by
 \beqnas
 G(X,Y)&=&\frac{1-XY}{(1-3XZ+3X^2\bZ-X^3)(1-3Y\bZ+3Y^2Z-Y^3)}
 \eeqnas
 \bpf
 This a a direct application of the computations provided in subsection~\ref{subsec.generating.gal} to this particular case.
 \epf
 \eprop

\section*{Aknowledgment }

I thank D. Bakry and N. Demni for many fruitful discussions and hints  about the content of this paper, together with a careful reading a  previous version of it.
 \bibliographystyle{amsplain}
\bibliography{bib.olfa}

\providecommand{\bysame}{\leavevmode\hbox to3em{\hrulefill}\thinspace}
\providecommand{\MR}{\relax\ifhmode\unskip\space\fi MR }
\providecommand{\MRhref}[2]{%
  \href{http://www.ams.org/mathscinet-getitem?mr=#1}{#2}
}
\providecommand{\href}[2]{#2}
\begin{thebibliography}{10}

\bibitem{bglbook}
D.~Bakry, I.~Gentil, and M.~Ledoux, \emph{Analysis and {G}eometry of {M}arkov
  {D}iffusion {O}perators}, Grund. Math. Wiss., vol. 348, Springer, Berlin,
  2013.

\bibitem{BOZ}
D.~Bakry, S.~Orevkov, and M.~Zani, \emph{Orthogonal polynomials and diffusions
  operators}, Preprint (2013).

\bibitem{beerends}
RJ~Beerends, \emph{Chebyshev polynomials in several variables and the radial
  part of the laplace-beltrami operator}, Transactions of the American
  Mathematical Society \textbf{328} (1991), no.~2, 779--814.

\bibitem{Bourb.Lie1}
Nicolas Bourbaki, \emph{Lie groups and {L}ie algebras. {C}hapters 1--3},
  Elements of Mathematics (Berlin), Springer-Verlag, Berlin, 1998, Translated
  from the French, Reprint of the 1989 English translation. \MR{1728312
  (2001g:17006)}

\bibitem{Bourb.Lie2}
\bysame, \emph{Lie groups and {L}ie algebras. {C}hapters 4--6}, Elements of
  Mathematics (Berlin), Springer-Verlag, Berlin, 2002, Translated from the 1968
  French original by Andrew Pressley. \MR{1890629 (2003a:17001)}

\bibitem{Bourb.lie3}
\bysame, \emph{Lie groups and {L}ie algebras. {C}hapters 7--9}, Elements of
  Mathematics (Berlin), Springer-Verlag, Berlin, 2005, Translated from the 1975
  and 1982 French originals by Andrew Pressley. \MR{2109105 (2005h:17001)}

\bibitem{dunn}
Ken~B Dunn and Rudolf Lidl, \emph{Generalizations of the classical chebyshev
  polynomials to polynomials in two variables}, Czechoslovak Mathematical
  Journal \textbf{32} (1982), no.~4, 516--528.

\bibitem{Caracter}
{{\"O}}mer E{\u{g}}ecio{\u{g}}lu, \emph{Algorithms for the character theory of
  the symmetric group}, E{UROCAL} '85, {V}ol.\ 2 ({L}inz, 1985), Lecture Notes
  in Comput. Sci., vol. 204, Springer, Berlin, 1985, pp.~206--224. \MR{826566}

\bibitem{Faraut}
Jacques Faraut, \emph{Analyse sur les groupes de lie: une introduction},
  Calvage \&amp; Mounet, 2006.

\bibitem{hum}
James~E Humphreys, \emph{Reflection groups and coxeter groups}, vol.~29,
  Cambridge university press, 1992.

\bibitem{James}
Gordon~Douglas James, Adalbert Kerber, PM~Cohn, and G~de~B Robinson, \emph{The
  representation theory of the symmetric group}, vol.~16, Cambridge University
  Press Cambridge, 1984.

\bibitem{Koorn1}
T.~Koornwinder, \emph{Orthogonal polynomials in two variables which are
  eigenfunctions of two algebraically independent partial differential
  operators. i.}, Nederl. Akad. Wetensch. Proc. Ser. A 77=Indag. Math.
  \textbf{36} (1974), 48--58.

\bibitem{Koorn2}
\bysame, \emph{Orthogonal polynomials in two variables which are eigenfunctions
  of two algebraically independent partial differential operators. ii.},
  Nederl. Akad. Wetensch. Proc. Ser. A 77=Indag. Math. \textbf{36} (1974),
  59--66.

\bibitem{Koorn}
\bysame, \emph{Two-variable analogues of the classical orthogonal polynomials},
  Theory and application of special functions, Proc. Advanced Sem., Math. Res.
  Center, Univ. Wisconsin, Madison, Wis., vol.~35, Academic Press, New--York,
  1975.

\bibitem{yvette}
Yvette Kosmann-Schwarzbach, \emph{Groupes et sym{\'e}tries: groupes finis,
  groupes et alg{\`e}bres de lie, repr{\'e}sentations}, Editions Ecole
  Polytechnique, 2005.

\bibitem{KrallS}
H.L. Krall and I.M. Sheffer, \emph{Orthogonal polynomials in two variables},
  Ann. Mat. Pura Appl. \textbf{76} (1967), 325--376.

\bibitem{T.Levy}
Thierry L{{\'e}}vy, \emph{Schur-{W}eyl duality and the heat kernel measure on
  the unitary group}, Adv. Math. \textbf{218} (2008), no.~2, 537--575.
  \MR{2407946 (2009g:15075)}

\bibitem{macdonald}
I.~G. Macdonald, \emph{Symmetric functions and orthogonal polynomials},
  University Lecture Series, vol.~12, American Mathematical Society,
  Providence, RI, 1998, Dean Jacqueline B. Lewis Memorial Lectures presented at
  Rutgers University, New Brunswick, NJ. \MR{1488699 (99f:05116)}

\bibitem{Mazet97}
O.~Mazet, \emph{Classification des semi-groupes de diffusion sur $\mathbb{R}$
  associ\'es \`a une famille de polyn\^omes orthogonaux}, S\'eminaire de
  probabilit\'es, Lectures notes in Mathematics, vol. 1655, Springer, 1997,
  pp.~40--54.

\bibitem{NSU}
A.F Nikiforov, S.K. Suslov, and V.B. Uvarov, \emph{Classical orthogonal
  polynomials of a discrete variable}, Springer Series in Computational
  Physics, Springer--Verlag, Berlin, 1991.

\bibitem{opdam}
E.M. Opdam, \emph{Root systems and hypergeometric functions, iii, iv},
  Compositio Math. \textbf{67} (1988), 21--49, 191--209.

\bibitem{Pinsky1}
Mark~A. Pinsky, \emph{{The Eigenvalues of an Equilateral Triangle}}, Siam
  Journal on Mathematical Analysis \textbf{11} (1980).

\bibitem{Pinsky}
Mark~A Pinsky, \emph{Completeness of the eigenfunctions of the equilateral
  triangle}, SIAM journal on mathematical analysis \textbf{16} (1985), no.~4,
  848--851.

\bibitem{Dunkel}
Jasper~V. Stokman, Charles~F. Dunkl, and Y.~Xu, \emph{{Orthogonal Polynomials
  of Several Variables}}, Journal of Approximation Theory \textbf{112} (2001),
  318--319.

\bibitem{Yoshida}
K~Yosida, \emph{Functional analysis. 1978}, 1978.

\end{thebibliography}
\end{document}